\documentclass[a4paper, 11pt]{amsart}

\usepackage{amsmath, tabularx}
\usepackage{amsthm}
\usepackage{amssymb}
\usepackage{amsfonts}
\usepackage{paralist}
\usepackage{aliascnt}
\usepackage{amscd}
\usepackage{blkarray}
\usepackage{mathbbol}
\usepackage[inner=3cm,outer=3cm, bottom=3.5cm]{geometry}
\usepackage{setspace}
\usepackage[colorlinks=true,linkcolor=blue,urlcolor=blue]{hyperref}
\usepackage[initials]{amsrefs}

\usepackage{tikz}
\usetikzlibrary{matrix}
\usetikzlibrary{arrows}

\BibSpec{collection.article}{%
	+{}  {\PrintAuthors}                {author}
	+{,} { \textit}                     {title}
	+{.} { }                            {part}
	+{:} { \textit}                     {subtitle}
	+{,} { \PrintContributions}         {contribution}
	+{,} { \PrintConference}            {conference}
	+{}  {\PrintBook}                   {book}
	+{,} { }                            {booktitle}
	+{,} { }                            {series}
	+{, vol.} { }                            {volume}
	+{,} { }                            {publisher}
	+{,} { \PrintDateB}                 {date}
	+{,} { pp.~}                        {pages}
	+{,} { }                            {status}
	+{,} { \PrintDOI}                   {doi}
	+{,} { available at \eprint}        {eprint}
	+{}  { \parenthesize}               {language}
	+{}  { \PrintTranslation}           {translation}
	+{;} { \PrintReprint}               {reprint}
	+{.} { }                            {note}
	+{.} {}                             {transition}
	+{}  {\SentenceSpace \PrintReviews} {review}
}

\newcommand{\kk}{\mathbb{k}}
\newcommand{\KK}{\mathbb{K}}
\newcommand{\bFF}{\mathbb{F}}

\newcommand{\CC}{\mathbb{C}}

\newcommand{\NN}{\normalfont\mathbb{N}}
\newcommand{\ZZ}{\normalfont\mathbb{Z}}

\newcommand{\MM}{{\normalfont\mathcal{M}}}

\newcommand{\mm}{{\normalfont\mathfrak{m}}}

\newcommand{\pp}{{\normalfont\mathfrak{p}}}

\newcommand{\Ker}{\normalfont\text{Ker}}

\newcommand{\Quot}{\normalfont\text{Quot}}

\newcommand{\Ann}{\normalfont\text{Ann}}

\newcommand{\Ass}{\normalfont\text{Ass}}

\newcommand{\Hom}{\normalfont\text{Hom}}

\newcommand{\gr}{\normalfont\text{gr}}
\newcommand{\Sol}{\normalfont\text{Sol}}

\newcommand{\Spec}{\normalfont\text{Spec}}
\newcommand{\Diff}{{\normalfont\text{Diff}}}
\newcommand{\DiffR}{\Diff_{R/A}}
\newcommand{\Diag}{\Delta_{R/A}}
\newcommand{\Princ}{P_{R/A}}

\newtheorem{theorem}{Theorem}[section]

\newtheorem{headthm}{Theorem}

\newaliascnt{headcor}{headthm}
\newtheorem{headcor}[headcor]{Corollary}
\aliascntresetthe{headcor}

\newaliascnt{headconj}{headthm}

\aliascntresetthe{headconj}

\newaliascnt{corollary}{theorem}
\newtheorem{corollary}[corollary]{Corollary}
\aliascntresetthe{corollary}

\newaliascnt{lemma}{theorem}
\newtheorem{lemma}[lemma]{Lemma}
\aliascntresetthe{lemma}

\newaliascnt{conjecture}{theorem}

\aliascntresetthe{conjecture}

\newaliascnt{proposition}{theorem}
\newtheorem{proposition}[proposition]{Proposition}
\aliascntresetthe{proposition}

\theoremstyle{definition}
\newaliascnt{definition}{theorem}
\newtheorem{definition}[definition]{Definition}
\aliascntresetthe{definition}

\newaliascnt{notation}{theorem}
\newtheorem{notation}[notation]{Notation}
\aliascntresetthe{notation}

\newaliascnt{example}{theorem}
\newtheorem{example}[example]{Example}
\aliascntresetthe{example}

\newaliascnt{examples}{theorem}

\aliascntresetthe{examples}

\newaliascnt{remark}{theorem}
\newtheorem{remark}[remark]{Remark}
\aliascntresetthe{remark}

\newaliascnt{problem}{theorem}

\aliascntresetthe{problem}

\newaliascnt{construction}{theorem}

\aliascntresetthe{construction}

\newaliascnt{setup}{theorem}
\newtheorem{setup}[setup]{Setup}
\aliascntresetthe{setup}

\newaliascnt{algorithm}{theorem}

\aliascntresetthe{algorithm}

\newaliascnt{observation}{theorem}

\aliascntresetthe{observation}

\newaliascnt{defprop}{theorem}

\aliascntresetthe{defprop}

\def\equationautorefname~#1\null{(#1)\null}
\def\sectionautorefname~#1\null{Section #1\null}
\def\subsectionautorefname~#1\null{\S #1\null}

\begin{document}

\title{Noetherian operators, primary submodules and symbolic powers}

\author{Yairon Cid-Ruiz}
\address{Max Planck Institute for Mathematics in the Sciences, Inselstra\ss e 22, 04103 Leipzig, Germany.}
\email{cidruiz@mis.mpg.de}
\urladdr{https://ycid.github.io}

\date{\today}

\begin{abstract}
	We give an algebraic and self-contained proof of the existence of the so-called \textit{Noetherian operators} for primary submodules over general classes of Noetherian commutative rings. 
	The existence of Noetherian operators accounts to provide an equivalent description of primary submodules in terms of differential operators. 
	As a consequence, we introduce a new notion of differential powers which coincides with symbolic powers in many interesting non-smooth settings, and so it could serve as a generalization of the Zariski-Nagata Theorem.
\end{abstract}

\subjclass[2010]{Primary: 13N10, 13N99, Secondary: 13A15.}

\keywords{differential operators, primary ideals, primary submodules, Noetherian operators, symbolic powers, differential powers.}

\maketitle

\section{Introduction}

The \textit{Fundamental Principle of Ehrenpreis and Palamodov} (see \cite{EHRENPREIS} and \cite{PALAMODOV}) is a celebrated theorem which states that all the solutions of a linear system of partial differential equations with constant coefficients can be represented as certain integrals of exponential-polynomial solutions.
Curiously enough, one of the main steps in the proof of this important theorem is to describe primary submodules of a finitely generated free module by using certain differential operators.
Following the terminology of Palamodov \cite{PALAMODOV}, these operators are commonly called \textit{Noetherian operators} in the literature.
In the case of polynomial rings over the complex numbers, one can find several proofs, that use algebraic and analytic techniques, for the existence of Noetherian operators, see, e.g., \cite[Chapter 8]{BJORK}, \cite{EHRENPREIS}, \cite{PALAMODOV}, \cite[\S 7.7]{HORMANDER}.

\medskip

It seems that a general characterization of primary submodules of a finitely generated free module over a polynomial ring with complex coefficients that used differential operators was obtained first in the form of Palamodov's Noetherian operators.
However, it turns out that Gr\"obner considered this problem before.
Gr\"obner described prime ideals and primary ideals having zero Krull dimension with the use of differential operators (see \cite{GROBNER_BOOK_AG_2, GROBNER_LIEGE, GROBNER_MATH_ANN}).
The terminology used by Palamodov (employing the term: \textit{Noetherian operators}; see \cite[Chapter IV, \S 3, page 161]{PALAMODOV}) is inspired by Noether's Fundamentalsatz (see, e.g., \cite[Chapter XIII, \S 96]{WAERDEN_ALG_II}).

 Subsequent algebraic approaches to characterize primary ideals and primary submodules with the use of differential operators were given in the following papers:

\begin{itemize}
	\item In \cite{BRUMFIEL_DIFF_PRIM}, for algebras over a field and under some assumptions, primary ideals were described by using differential operators from the algebra to the residue field of the corresponding prime.
	In \cite{BRUMFIEL_DIFF_PRIM}, there is also a description of primary submodules.
	\item In \cite{BOMMER}, for algebras essentially of finite type over a perfect field, primary ideals to regular prime ideals were described by using derivations.
	\item In \cite{OBERST_NOETH_OPS}, the existence of Noetherian operators was proven for the case of polynomial rings over any field.
\end{itemize} 
But, in \cite{BRUMFIEL_DIFF_PRIM} and \cite{BOMMER} there is no reference to the previously obtained results of \cite{EHRENPREIS} and \cite{PALAMODOV}.

\medskip

The main purpose of this paper is to provide an algebraic and self-contained development of the existence of Noetherian operators for general classes of Noetherian commutative rings.
 In other words, for large classes of Noetherian commutative rings, we show the existence of an equivalent notion of primary submodules that depends upon differential operators.
As a consequence, we introduce a new notion of differential powers which could serve as a generalization of the Zariski-Nagata Theorem (\cite{ZARISKI,NAGATA_LOCAL_RINGS}) for non-smooth settings.
Unless specified otherwise, in this paper all rings are commutative.
  
\medskip 
  
Next is a summary of the  main results of this paper, to simplify the exposition here in the introduction, below we state them only for the case of primary ideals.
For the rest of the introduction, let $A$ be a Noetherian integral domain and $R$ be a Noetherian ring such that $A \subset R$. 
Denote by $\Quot(A)$ the field of fractions of $A$.

\medskip

The first main result deals with the problem of describing primary ideals as solution sets of certain differential operators (see \autoref{def_diff_ops}).

\begin{headthm}[{\autoref{thm_Pimary_ideals_as_Sol}, \autoref{thm_Pimary_submod_as_Sol}}]
	\label{thmA}
	Let $R$ be a Noetherian ring and $A$ be a subring, such that $A$ is a Noetherian integral domain.
	Let $\pp \in \Spec(R)$ be a prime ideal in $R$ such that $\pp \cap A = 0$.
	Then, the following statements hold:
	\begin{enumerate}[(i)]
		\item Suppose that $R_\pp/\pp R_\pp \otimes_{\Quot(A)} R_\pp$ is a Noetherian ring.
		If $I \subset R$ is a $\pp$-primary ideal in $R$, then there exists an $(R \otimes_A R)$-submodule $\mathcal{E} \subseteq \DiffR\left(R,R_\pp/\pp R_\pp\right)$ such that 
		$$
		I = \big\lbrace f \in R \mid \delta(f) = 0 \text{ for all } \delta \in \mathcal{E} \big\rbrace
		$$
		and $\mathcal{E}$ is finitely generated via its natural left structure as an $R_\pp$-module.
		\item Suppose that $R$ is essentially of finite type over $A$ and $N$ is a finitely generated torsion-free module over $R/\pp$.
		\begin{enumerate}[(a)]
			\item If $I \subset R$ is a $\pp$-primary ideal in $R$, then there exists an $(R \otimes_A R)$-submodule $\mathcal{E} \subseteq \DiffR(R,N)$ such that 
			$$
			I = \big\lbrace f \in R \mid \delta(f) = 0 \text{ for all } \delta \in \mathcal{E} \big\rbrace
			$$
			and $\mathcal{E}$ is finitely generated via its natural left structure as an $R$-module.
			\item If $\Quot(A) \hookrightarrow R_\pp/\pp R_\pp$ is a separable field extension, which holds whenever $\Quot(A)$ is perfect, then for any $\pp$-primary ideal $I \subset R$ containing $\pp^{n+1}$, there exists an $(R \otimes_A R)$-submodule $\mathcal{E} \subseteq \DiffR^n(R,N)$ such that 
			$$
			I = \big\lbrace f \in R \mid \delta(f) = 0 \text{ for all } \delta \in \mathcal{E} \big\rbrace.
			$$
		\end{enumerate}				
	\end{enumerate}
\end{headthm}

We now discuss the concept of Noetherian operators in the context of \autoref{thmA}.
In \autoref{thmA}$(i)$, if $\{ \delta_1,\ldots,\delta_q \}$ is a set of generators of $\mathcal{E}$ as a left $R_\pp$-module, we obtain that $I = \big\lbrace f \in R \mid \delta_i(f) = 0 \text{ for all } 1 \le i \le q \big\rbrace$ (see \autoref{rem:finite_gen_Noeth_ops}).
Accordingly, we shall say that $\{ \delta_1,\ldots,\delta_q \}$ is a set of Noetherian operators for the $\pp$-primary ideal $I$.
Similarly, in \autoref{thmA}$(ii)$, a finite set of generators of $\mathcal{E}$ as a left $R$-module will be called a set of Noetherian operators for the $\pp$-primary ideal $I$.

\medskip

The above \autoref{thmA} was initially inspired by the important results of \cite{BRUMFIEL_DIFF_PRIM}.
But, we generalize the main results of \cite{BRUMFIEL_DIFF_PRIM} in two ways: we do not assume that $A$ is a field and we use more general types differential operators (not just differential operators  in $\DiffR(R, R_\pp/\pp R_\pp)$).

\medskip

If we assume that $A$ is a field and $R$ is a polynomial ring over $A$, from \cite{OBERST_NOETH_OPS} we know that we can use differential operators in $\DiffR(R,R)$ to describe primary ideals.
Since there is much more literature and interest in the  differential operators in $\DiffR(R,R)$ (for instance, when $A$ is a field of characteristic zero and $R$ is a polynomial ring over $A$,  $\DiffR(R,R)$ is referred to as the Weyl algebra), it is natural to ask when the results of \autoref{thmA} can be stated by using differential operators in $\DiffR(R,R)$. 
Our second main result shows that such a description is possible under certain smooth settings.
We point out that this statement may not hold without any assumption of smoothness (see \autoref{exam_cubic}).

\begin{headcor}[{\autoref{cor_prim_ideals_sol_Diff_R_R}, \autoref{cor_prim_submods_sol_Diff_R_R}}]
	\label{corB}
	Let $A$ be a Noetherian integral domain and $R$ be an $A$-algebra formally smooth and essentially of finite type over $A$ such that $A \subset R$.
	Let $\pp \in \Spec(R)$ be a prime ideal in $R$ such that $\pp \cap A = 0$.
	Then, the following statements hold:
	\begin{enumerate}[(i)]
		\item If $I \subset R$ is a $\pp$-primary ideal in $R$, then there exists an $(R \otimes_A R)$-submodule $\mathcal{E} \subseteq \DiffR(R,R)$ such that 
		$$
		I = \big\lbrace f \in R \mid \delta(f) \in \pp \text{ for all } \delta \in \mathcal{E} \big\rbrace
		$$
		and $\mathcal{E}$ is finitely generated via its natural left structure as an $R$-module.
		\item If $\Quot(A) \hookrightarrow R_\pp/\pp R_\pp$ is a separable field extension, which holds whenever $\Quot(A)$ is perfect, then for any $\pp$-primary ideal $I \subset R$ containing $\pp^{n+1}$, there exists an  $(R \otimes_A R)$-submodule $\mathcal{E} \subseteq \DiffR^n(R,R)$ such that 
		$$
		I = \big\lbrace f \in R \mid \delta(f) \in \pp \text{ for all } \delta \in \mathcal{E} \big\rbrace.
		$$
	\end{enumerate}
\end{headcor}

In the context of \autoref{corB}, the Noetherian operators are given by using Palamodov's original terminology;  if $\{ \delta_1,\ldots,\delta_q \}$ is a set of generators of $\mathcal{E}$ as a left $R$-module, we obtain that $I = \big\lbrace f \in R \mid \delta_i(f) \in \pp \text{ for all } 1 \le i \le q \big\rbrace$.
Here we shall also say that $\{ \delta_1,\ldots,\delta_q\}$ is a set of Noetherian operators for the $\pp$-primary ideal $I$.

\medskip

The symbolic powers of a prime ideal are very special primary ideals that have received a lot of attention in the areas of Algebraic Geometry and Commutative Algebra (see \cite{RESUME_SYMB_HUNEKE_ET_AL}). 
The Zariski-Nagata Theorem (see, e.g., \cite{ZARISKI}, \cite{NAGATA_LOCAL_RINGS}, \cite[Theorem 3.14]{EISEN_COMM}, \cite{RESUME_SYMB_HUNEKE_ET_AL}) is a fundamental result that, in the case where $R$ is a polynomial ring over a perfect field, describes the $n$-th symbolic power of a given prime ideal as the polynomials that vanish to order greater than or equal to $n$ on the corresponding variety.
To extend this study in other rings, for any ideal $I \subset R$, one has the following ideals
$$
I^{{\langle n \rangle}_A} = \big\lbrace f \in R \mid \delta(f) \in I \text{ for all } \delta \in \DiffR^{n-1}(R,R)  \big\rbrace,
$$
dubbed as differential powers (see \cite{RESUME_SYMB_HUNEKE_ET_AL}). 
These differential powers have sparked attention and renewed interest in extending the Zariski-Nagata Theorem (see, e.g., \cite{brenner2018quantifying, ZAR_NAG_DESTEFANI_ET_AL}).

In this paper, we propose a new notion of differential powers which seems to be better suited to describe symbolic powers, especially because it coincides with symbolic powers in many interesting non-smooth settings.
From \cite{EISENBUD_HOCHSTER}, it has been long known that some assumption of smoothness is needed to extend the original statement of the Zariski-Nagata Theorem.
Also, see \autoref{exam_cubic}.
Recently, in \cite[Theorem 10.2]{brenner2018quantifying}, it has been shown that the conclusion of the original Zariski-Nagata Theorem characterizes smoothness. 

For any ideal $I \subset R$, we introduce the following new notion of differential powers
$$
I^{{\lbrace n \rbrace}_A} = \big\lbrace f \in R \mid \delta(f) = 0 \text{ for all } \delta \in \DiffR^{n-1}(R,R/I)  \big\rbrace.
$$
In our last main result, we relate the two above notions of differential powers with symbolic powers.

\begin{headthm}[{\autoref{thm_symb_powers_diff_powers}}]
	\label{thmC}
	Let $A$ be a Noetherian integral domain and $R$ be an $A$-algebra essentially of finite type over $A$ such that $A \subset R$.
	Let $\pp \in \Spec(R)$ be a prime ideal in $R$ such that $\pp \cap A = 0$.
	Then, the following statements hold:
	\begin{enumerate}[(i)]
		\item $\pp^{(n)} \subseteq \pp^{{\{n\}}_A} \subseteq \pp^{{\langle n \rangle}_A}$.
		\item If $\Quot(A) \hookrightarrow R_\pp/\pp R_\pp$ is a separable field extension, which holds whenever $\Quot(A)$ is perfect, then
		$$
		\pp^{(n)} = \pp^{{\{n\}}_A}.
		$$ 
		\item If $R$ is formally smooth over $A$, then 
		$$
		\pp^{{\{n\}}_A} = \pp^{{\langle n \rangle}_A}.
		$$		
	\end{enumerate}
\end{headthm}

The basic outline of this paper is as follows.
In \autoref{section_diff_ops}, we recall some basic results on differential operators that will be needed throughout the rest of the paper.
In \autoref{sect_primary_submodules}, we prove \autoref{thmA} and \autoref{corB}.
In \autoref{sect_Zar_Nag}, we prove \autoref{thmC}.
In \autoref{sect_examples}, we provide some examples and computations.

\section{Differential operators}
\label{section_diff_ops}

During this short section we recall some basic notions regarding differential operators.
A general and complete reference on the topic of differential operators is \cite[\S 16]{EGAIV_IV}.

Throughout this paper, unless specified otherwise,  a ring is always assumed to be a commutative ring.
The following setup will be used during the present section.

\begin{setup}
	Let $R$ be a ring and $A$ be a subring.
\end{setup}

For two $R$-modules $M$ and $N$, we regard $\Hom_A(M, N)$ as an $(R\otimes_A R)$-module, by setting 
$$
\left((r \otimes_A s) \delta\right)(m) = r \delta(sm) \quad \text{ for all } \delta \in \Hom_A(M, N), \; m \in M,\; r,s \in R. 
$$ 
We use the bracket notation $[\delta,r](m) = \delta(rm)-r\delta(m)$ for all $\delta \in \Hom_A(M, N)$, $r \in R$ and $m \in M$.
The $A$-linear differential operators form an $(R \otimes_A R)$-submodule of $\Hom_A(M, N)$ and are defined inductively as follows.

\begin{definition}
	\label{def_diff_ops}
	Let $R$ be a  ring and $A$ be a subring.
	Let $M, N$ be $R$-modules.
	The $n$-th order $A$-linear differential operators $\DiffR^n(M, N) \subseteq \Hom_A(M, N)$ from $M$ to $N$ are defined inductively by:
	\begin{enumerate}[(i)]
		\item $\DiffR^{-1}(M,N) := 0$.
		\item $\DiffR^{n}(M, N) := \big\lbrace \delta \in \Hom_A(M,N) \mid [\delta, r] \in \DiffR^{n-1}(M, N) \text{ for all } r \in R \big\rbrace$.
	\end{enumerate}
	The $A$-linear differential operators from $M$ to $N$ are given by 
	$$
	\DiffR(M, N) := \bigcup_{n=0}^\infty \DiffR^n(M,N).
	$$
\end{definition}

Right from the definition we obtain that $\DiffR^0(M, N)=\Hom_R(M, N)$.
We endow $\DiffR^n(M, N)$ with a structure of $R$-module given by the natural map $R \rightarrow R \otimes_A R$, $r \mapsto r \otimes_A 1$ and the $(R\otimes_A R)$-module structure of $\DiffR^n(M, N)$ (see \autoref{prop_Diff_ann_diag} below), that is, $R$ acts as post-composition of maps over $\DiffR^n(M, N)$: for $r \in R, \delta \in \DiffR^n(M, N)$, $r\delta$ is the differential operator given by $(r\delta)(m) = r \delta(m)$ for all $m \in M$.

To describe differential operators, a fundamental idea is to study the module of principal parts.
Consider the multiplication map 
\begin{equation*}
	\mu : R \otimes_A R \rightarrow R, \quad 	r \otimes_A s \mapsto rs,
\end{equation*}
and denote the kernel of this map as the ideal $\Diag \subset R \otimes_A R$. 
By making a simple induction argument, one has the following equivalent description.

\begin{proposition}[{\cite[Proposition 2.2.3]{AFFINE_HOPF_I}}]
	\label{prop_Diff_ann_diag}
	Let $M, N$ be $R$-modules.
	Then, $\DiffR^n(M, N)$ is the $(R \otimes_A R)$-submodule of $\Hom_A(M, N)$ annihilated by $\Diag^{n+1}$.
\end{proposition}

\begin{definition}
	Let $M$ be an $R$-module.
	The module of $n$-th principal parts is the $(R \otimes_A R)$-module defined by 
	$$
	\Princ^n(M) := \frac{R \otimes_A M}{\Diag^{n+1}  \left(R \otimes_A M\right)}.
	$$
	For simplicity of notation, $\Princ^n(R)$ is denoted as $\Princ^n$.
\end{definition}

\begin{remark}
	\label{rem_isom_mod_princ_parts}
	We have an isomorphism of $(R \otimes_A R)$-modules 
	$$
	\Princ^n(M) \,=\, \frac{R \otimes_A M}{\Diag^{n+1}  \left(R \otimes_A M\right)} \,\cong\, \frac{R \otimes_A R}{\Diag^{n+1}} \otimes_R M \,=\, \Princ^n \otimes_R M
	$$
	where the tensor product is taken via the right $R$-module structure of $\frac{R \otimes_A R}{\Diag^{n+1}}$.
\end{remark}

Unless we specify otherwise, whenever we consider $\Princ^n(M)$ as an $R$-module, we do it so by setting that $R$ acts over the left factor of $R \otimes_A M$, that is, for $r,s \in R$ and $m \in M$ we have
$$
r \big(\overline{s \otimes_A m}\big) = \overline{rs \otimes_A m},
$$
where $\overline{s \otimes_A m}$ represents the residue class of $s \otimes_A m \in R \otimes_A M$ in $\Princ^n(M)$.

It turns out that $\DiffR^n$ can be seen as a representable functor, as follows. 
If there are $R$-homomorphisms $f: M^\prime \rightarrow M$ and $g : N \rightarrow N^\prime$, then the $A$-homomorphism 
$$
\Hom_A(f, g): \Hom_A(M, N) \rightarrow \Hom_A(M^\prime, N^\prime), \quad \delta \mapsto g \circ \delta \circ f 
$$
is naturally an $(R \otimes_A R)$-homomorphism which is compatible with the bracket operation, that is, $\left[g \circ \delta \circ f, r\right] = g \circ \left[\delta, r\right] \circ f$ for all $r \in R$.
So, from \autoref{prop_Diff_ann_diag}, we have that $\DiffR^n(\bullet, \bullet)$ is a functor
\begin{equation}
	\label{eq_Diff_functor}
	\DiffR^n(f, g): \DiffR^n(M, N) \rightarrow \DiffR^n(M^\prime, N^\prime), \quad \delta \mapsto g \circ \delta \circ f. 
\end{equation}
The adjointness of Hom and tensor gives an isomorphism 
$$
\Hom_R\left(R \otimes_A M, N\right) \cong \Hom_A(M, N)
$$
of $(R \otimes_A R)$-modules, where the $\left(R \otimes_A R\right)$-module structure of $\Hom_R\left(R \otimes_A M, N\right)$ is given by setting 
$$
\left(\left(r \otimes_A s\right) \psi\right)\left(t \otimes_A m\right) = \psi\left(rt \otimes_A sm\right) = rt \psi\left(1 \otimes_A sm\right)
$$
for all $\psi \in \Hom_R\left(R \otimes_A M, N\right)$, $r,s,t \in R$, $m \in M$.
Indeed, one has that for any $\psi \in \Hom_R\left(R\otimes_A M, N\right)$ one can define $\varphi \in \Hom_A(M, N)$ as $\varphi(m)=\psi(1 \otimes_A m)$, and, in the other direction, for any $\varphi \in \Hom_A(M,N)$ one can define $\psi \in \Hom_R\left(R \otimes_A M, N\right)$  as $\psi(r \otimes_A m) = r\varphi(m)$; also, these maps are inverse to each other.
Therefore, by using the universal map 
$$
d^n : M \rightarrow \Princ^n(M), \qquad m \in M \mapsto \overline{1 \otimes_A m} \in \Princ^n(M),
$$
and \autoref{prop_Diff_ann_diag}, we can obtain the following description of differential operators.

\begin{proposition}[{\cite[Proposition 16.8.4]{EGAIV_IV}, \cite[Theorem 2.2.6]{AFFINE_HOPF_I}}]
	\label{prop_represen_diff_opp}
	Let $n\ge 0$ and $M, N$ be $R$-modules.
	Then, the map 
	\begin{align*}
		{\left(d^n\right)}^* : \Hom_R\left(\Princ^n(M), N\right) &\xrightarrow{\cong} \DiffR^n(M, N) \\
		\varphi &\mapsto \varphi \circ d^n
	\end{align*}
	induces an isomorphism of $R$-modules.
\end{proposition}

Finally, the following lemma gathers some general properties of differential operators.
Most of these basic results have appeared elsewhere (see, e.g.,  \cite{brenner2018quantifying, KSMITH_VAN_DER_BERGH}), but one includes them for the sake of completeness.

\begin{lemma}
	\label{lem_gen_props_diff_ops}
	Suppose that we have an inclusion of rings $A \subseteq B \subseteq R$, and that $W \subset R$ is a multiplicatively closed subset in $R$.
	Let $M, N$ be $R$-modules.
	Then, for $n \ge 0$, the following statements hold:
	\begin{enumerate}[(i)]
		\item $\DiffR^n(M, \bullet)$ is a left exact covariant functor and $\DiffR^n(\bullet, N)$ is a left exact contravariant functor.
		\item $\Diff_{R/B}^n(M, N) \subseteq \DiffR^n(M, N)$.
		\item If $N$ is a module over $W^{-1}R$, then  
		$$
		\DiffR^n(M, N) \cong \Diff_{W^{-1}R/A}^n\left(W^{-1}M, N\right) \cong \Diff_{W^{-1}R/{(W \cap A)}^{-1}A}^n\left(W^{-1}M, N\right)
		$$
		are isomorphisms of  $R$-modules.
		\item If $\Princ^n(M)$ is finitely presented as an $R$-module, then 
		\begin{align*}
			W^{-1}\DiffR^n(M, N) &\cong \Diff_{W^{-1}R/A}^n\left(W^{-1}M, W^{-1}N\right)\\
			 &\cong \Diff_{W^{-1}R/{(W \cap A)}^{-1}A}^n\left(W^{-1}M, W^{-1}N\right)
		\end{align*}
		are isomorphisms of $R$-modules.
		\item If ${\left(M_\lambda\right)}_{\lambda \in \Lambda}$ is a family of $R$-modules, then 
		$$
		\DiffR^n\left(\bigoplus_{\lambda \in \Lambda}M_\lambda, N\right) \;\cong\;  \prod_{\lambda \in \Lambda} \DiffR^n\left(M_\lambda, N\right)
		$$
		is an isomorphism of $R$-modules.
		\item If ${\left(N_\lambda\right)}_{\lambda \in \Lambda}$ is a family of $R$-modules, then 
		$$
		\DiffR^n\left(M, \prod_{\lambda \in \Lambda}N_\lambda\right) \;\cong\;  \prod_{\lambda \in \Lambda} \DiffR^n\left(M, N_\lambda\right)
		$$
		is an isomorphism of $R$-modules.
	\end{enumerate}
	\begin{proof}		
		$(i)$ 
		Here we only prove the statement about $\DiffR^n(\bullet, N)$, the proof of the other claim about $\DiffR^n(M,\bullet)$ is completely similar.
		
		Suppose we have an exact sequence $M^\prime \rightarrow M \rightarrow M^{\prime\prime} \rightarrow 0$ of $R$-modules, then the right-exactness of tensor product and \autoref{rem_isom_mod_princ_parts} yield the exact sequence 
		$$
		\Princ^n\left(M^\prime\right) \rightarrow \Princ^n\left(M\right) \rightarrow \Princ^n\left(M^{\prime\prime}\right) \rightarrow 0.
		$$
		The left-exactness of $\Hom_R(\bullet, N)$ gives the exact sequence 
		$$
		0 \rightarrow \Hom_R\left(\Princ^n\left(M^{\prime\prime}\right), N\right) \rightarrow \Hom_R\left(\Princ^n\left(M\right), N\right) \rightarrow \Hom_R\left(\Princ^n\left(M^{\prime}\right), N\right),
		$$
		and so the result follows from \autoref{prop_represen_diff_opp}.
		
		$(ii)$ Note that there is a canonical surjection 
		$$
		\Princ^n(M) = \frac{R \otimes_A M}{\Diag^{n+1}  \left(R \otimes_A M\right)} \,\twoheadrightarrow\, P_{R/B}^n(M) = \frac{R \otimes_B M}{\Delta_{R/B}^{n+1}  \left(R \otimes_B M\right)}.
		$$
		Hence, \autoref{prop_represen_diff_opp} gives the inclusion $\Diff_{R/B}^n(M, N) \subseteq \DiffR^n(M, N)$.

		$(iii)$
		By \autoref{prop_represen_diff_opp} and the Hom-tensor adjunction, we get the isomorphisms 
		\begin{align}
		\label{eq_diff_op_N_is_localized}
		\begin{split}
		\DiffR^n(M, N) &\cong \Hom_R\left(\Princ^n(M), N\right) \\
		&\cong \Hom_R\left(\Princ^n(M), \Hom_{W^{-1}R}\left(W^{-1}R, N\right)\right)\\
		&\cong \Hom_{W^{-1}R}\left(W^{-1}R \otimes_R \Princ^n(M), N\right).	
		\end{split}			
		\end{align}
		From \cite[Proposition 2.15]{brenner2018quantifying} (\cite[Proposition 16.4.14]{EGAIV_IV}) we have the isomorphisms 
		\begin{equation}
		\label{eq_localization_Prin}
		W^{-1}R \otimes_R \Princ^n \cong P_{W^{-1}R/ A}^n \cong P_{W^{-1}R/ {(W\cap A)}^{-1}A}^n.
		\end{equation}
		By summing up \autoref{eq_diff_op_N_is_localized}, \autoref{eq_localization_Prin}, \autoref{rem_isom_mod_princ_parts} and \autoref{prop_represen_diff_opp}, the result is obtained.
		
		$(iv)$	It follows from \autoref{rem_isom_mod_princ_parts},  \autoref{prop_represen_diff_opp} and \autoref{eq_localization_Prin} because $\Hom$ and localization commute 
		$$
		W^{-1}R \otimes_R \Hom_R\left(\Princ^n(M), N\right) \cong \Hom_{W^{-1}R}\left(W^{-1}R \otimes_R \Princ^n(M), W^{-1}N\right)
		$$
		when $\Princ^n(M)$ is a finitely presented $R$-module (see, e.g., \cite[Theorem 7.11]{MATSUMURA}).

		$(v)$ Since tensor products commute with direct sums, from \autoref{rem_isom_mod_princ_parts} we obtain the isomorphism
		$$
		\Princ^n\left(\bigoplus_{\lambda \in \Lambda}M_\lambda\right) \cong \bigoplus_{\lambda \in \Lambda} \Princ^n\left(M_\lambda\right).
		$$
		Therefore, from the isomorphism 
		$$
		\Hom_R\left(\bigoplus_{\lambda \in \Lambda}\Princ^n\left(M_\lambda\right), N\right) \cong \prod_{\lambda \in \Lambda} \Hom_R\left(\Princ^n\left(M_\lambda\right), N\right)
		 $$
		 and \autoref{prop_represen_diff_opp}, the result follows.
		
		$(vi)$ It follows similarly to $(v)$.
	\end{proof}
\end{lemma}

\section{Primary submodules as the zero set of differential operators}
\label{sect_primary_submodules}

In this section we show an algebraic and more general version of the existence of Noetherian operators. 
In a quite unrestrictive setting, we shall prove that any primary submodule of a finitely generated module can be obtained as the solution set of certain differential operators.

The following setup is used throughout this section. 

\begin{setup}
	\label{setup_primary_submod}
	Let $R$ be a Noetherian ring and $A$ be a subring, such that $A$ is a Noetherian integral domain.
\end{setup}

Note that the above setup is more general than the ones of \cite{BRUMFIEL_DIFF_PRIM} and \cite{OBERST_NOETH_OPS} because in both of those papers it is assumed that $A$ is a field.
Also, we do not assume any finiteness condition of $R$ over $A$, which may lead to cases where $R \otimes_A R$ is not a Noetherian ring and the modules of principal parts $\Princ^n$ are not finitely generated modules over $R$.

For simplicity we use the notation below.

\begin{notation}
	(i)	 Let $M, N$ be $R$-modules. 
		For a subset $\mathcal{E} \subseteq \DiffR(M, N)$, we set 
		$$
		\Sol(\mathcal{E}) \,:=\, \big\lbrace m \in M \mid \delta(m) = 0 \text{ for all } \delta \in \mathcal{E} \big\rbrace 
		\,=\, \bigcap_{\delta \in \mathcal{E} } \Ker(\delta).
		$$
		
	\noindent	
	(ii) Denote by $\Quot(A)$ the field of fractions of $A$.
	
	\smallskip
	\noindent		
	(iii) For any $\pp \in \Spec(R)$, we denote by $k(\pp)$ the residue field 
		$
		k(\pp):=R_\pp/\pp R_\pp = \Quot(R/\pp)
		$
		at $\pp$.		
\end{notation}

To obtain finer results we will need to assume that $R$ is essentially of finite type over $A$, i.e., $R$ is the localization of a finitely generated $A$-algebra.
The remark below contains some easy consequences of that additional assumption.

\begin{remark}
	\label{rem_conseq_essent_finite_type}
	Additionally, if $R$ is essentially of finite type over $A$, then the following statements hold:
	\begin{enumerate}[(i)]
		\item $k(\pp)$ is a finitely generated field over $\Quot(A)$.
		\item $R \otimes_A R$ is a Noetherian ring.
		\item For all $n \ge 0$, $\Princ^n$ is a finitely generated $R$-module, where, as before, the $R$-module structure comes from the left factor of $R \otimes_A R$. 
	\end{enumerate}
	\begin{proof}
		$(i)$ This is quite clear. 
		
		$(ii)$ Let $T$ be a finitely generated $A$-algebra such that $R$ is a localization of $T$.
		Then, $T \otimes_A T$ is Noetherian, and $R \otimes_A R$, being a localization of $T \otimes_A T$, is also Noetherian.
		
		$(iii)$ From part $(ii)$, each ideal $\Diag^n \subset R \otimes_A R$ is finitely generated, and so it follows that $\Diag^n/\Diag^{n+1}$ is a finitely generated module over $ \left(R \otimes_A R\right) / \Diag \xrightarrow{\cong} R$.
		Therefore, the following exact sequence
		$$
		0 \rightarrow \Diag^n/\Diag^{n+1} \rightarrow \Princ^{n} \rightarrow \Princ^{n-1} \rightarrow 0
		$$
		and an induction argument imply the result.
	\end{proof}
\end{remark}

For $R$-modules $M$ and $N$, the differential operators $\DiffR(M, N)$ have the filtration 
$$
\DiffR^0(M, N) \subseteq \cdots \subseteq \DiffR^{n-1}(M, N) \subseteq \DiffR^n(M, N) \subseteq \cdots, 
$$
which is of utmost importance. 
In particular, we will focus on the fact that $(R \otimes_A R)$-submodules of $\DiffR^n(M, N)$  are stable with respect to this filtration in the following easy sense. 

\begin{remark}
	For an $(R \otimes_A R)$-submodule $\mathcal{E}$ of $\DiffR^n(M, N)$ we have that 
	$$
	\left[\delta, r\right] = \left( 1 \otimes_A r -  r \otimes_A 1\right) \delta  \;\in\; \mathcal{E} \cap \DiffR^{k-1}(M, N)
	$$ for all $\delta \in \mathcal{E} \cap \DiffR^{k}(M, N)$, $r \in R$ and $k \le n$.
\end{remark}

The proposition below contains some basic properties of  $(R\otimes_A R)$-submodules of $\DiffR^n(M, N)$.
This proposition can also be found in \cite[Proposition 1.3]{BRUMFIEL_DIFF_PRIM}, but we include a proof for the sake of completeness and because in the proof of \cite[Proposition 1.3]{BRUMFIEL_DIFF_PRIM} there is a small gap (in the step of showing that $\pp^{n+1}M \subseteq \Sol(\mathcal{E})$ below).

\begin{proposition}[{\cite[Proposition 1.3]{BRUMFIEL_DIFF_PRIM}}]
	\label{prop_sol_closed_submods_diff}
	Let $M,N$ be $R$-modules and $\mathcal{E} \subseteq \DiffR^n(M, N)$ be an $(R \otimes_A R)$-submodule.
	Then, the following statements hold:
	\begin{enumerate}[(i)]
		\item If $\mathcal{E} \neq 0$, then $\mathcal{E} \, \cap\, \Hom_R(M, N) \neq 0$.
		\item $\Sol(\mathcal{E})$ is an $R$-submodule of $M$.
		Additionally, if $\mathcal{E} \neq 0$, then $\Sol(\mathcal{E}) \subsetneq M$.
		\item Let $\pp \in \Spec(R)$ be a prime ideal in $R$ and suppose that $N$ is a torsion-free module over $R/\pp$.
		Then, $\Sol(\mathcal{E})$ is a $\pp$-primary $R$-submodule of $M$ and $\pp^{n+1}M \subseteq  \Sol(\mathcal{E})$.
	\end{enumerate}
	\begin{proof}
		$(i)$ Suppose that $i$ is the least integer such that $\mathcal{E} \cap \DiffR^i(M, N) \neq 0$ and choose any $0 \neq \delta \in \mathcal{E} \cap \DiffR^i(M, N)$.
		Thus, for all $r \in R$, we have that $\left[\delta, r\right]=0$, which implies that $\delta \in \Hom_R(M, N)$,  and so $i =0$.
		
		$(ii)$	For all $r \in R, m \in M, \delta \in \mathcal{E}$, we have the equation
		$$
		\delta(rm)=r\delta(m) + \left[\delta, r\right](m).
		$$ 
		Hence, $m \in \Sol(\mathcal{E})$ implies that $rm \in \Sol(\mathcal{E})$, and so $\Sol(\mathcal{E})$ is an $R$-submodule.
		
	 	If $\mathcal{E} \neq 0$, then there exists $0 \neq \delta \in  \mathcal{E} \subset \Hom_A(M, N)$, which yields $\Ker(\delta) \subsetneq M$ and so $\Sol(\mathcal{E}) \subsetneq M$.
	
		$(iii)$ First, we show by induction on $n$ that $\pp^{n+1}M \subseteq  \Sol(\mathcal{E})$.
		For any $\delta \in \Hom_R(M, N)$, since $N$ is an $R/\pp$-module, we obtain that $\delta(\pp M) = \pp\delta(M) = 0$, and so $\pp M \subseteq \Sol(\mathcal{E})$ whenever $\mathcal{E} \subseteq \DiffR^0(M, N) = \Hom_R(M, N)$. 
		
		Now, suppose that $n > 0$.
		For $r \in \pp, s \in \pp^n, m \in M, \delta \in \DiffR^n(M, N)$, we obtain the equation
		$$
		\delta(rsm) = r\delta(sm) + \left[\delta, r\right](sm) = \left[\delta, r\right](sm).
		$$
		Since $\left[\delta, r\right]\in  \DiffR^{n-1}(M,N)$, the induction step yields $\delta(rsm) = \left[\delta, r\right](sm)=0$.
		Therefore, if $\mathcal{E} \subseteq \DiffR^{n}(M,N)$, then it follows that $\pp^{n+1}M \subseteq \Sol(\mathcal{E})$.
		
		Finally, we need to show that, if $rm \in \Sol(\mathcal{E})$ and $m \not\in \Sol(\mathcal{E})$, then $r \in \pp$.
		If $m \not\in \Sol(\mathcal{E})$, let $\delta \in \mathcal{E} \cap \DiffR^i(M, N)$ be a differential operator with smallest possible order $i$ that satisfies $\delta(m) \neq 0$.
		Thus, we have the equation $0 = \delta(rm)=r\delta(m) + \left[\delta,r\right](m)=r\delta(m)$. 
		The fact that $N$ is torsion-free over $R/\pp$ implies that $r \in \pp$, and so we get that $\Sol(\mathcal{E})$ is a $\pp$-primary $R$-submodule.
	\end{proof}
\end{proposition}

Now, before proceeding to the proof of the main theorem, which provides a converse for \autoref{prop_sol_closed_submods_diff}, there are some steps to reduce the problem to simpler situations.

The following proposition will be used when $R$ is essentially of finite type over $A$.

\begin{proposition}
	\label{prop_reduction_to_zero_dim}
	Assume \autoref{setup_primary_submod} with the additional condition that $R$ is essentially of finite type over $A$.
	Let $M$ be a finitely generated $R$-module.
	Let $\pp \in \Spec(R)$ be a prime ideal in $R$ such that $\pp \cap A = 0$, and $N$ be a torsion-free $R/\pp$-module.
	Then, the following statements hold:
	\begin{enumerate}[(i)]		
		\item If $\Quot(A) \hookrightarrow k(\pp)$ is separable, then there exists a field $\KK$ such that $\Quot(A) \subseteq \KK \subseteq R_\pp$ and $\KK \hookrightarrow k(\pp)$ is a separable finite field extension.
		
		\item Let $\KK$ be a field such that $\Quot(A) \subseteq \KK \subseteq R_\pp$.
		Then, for any $\left(R_\pp\otimes_\KK R_\pp\right)$-submodule $\mathcal{E}^\prime \subseteq \Diff_{R_\pp/\KK}^n\left(M_\pp, N_\pp\right)$, there exists an $(R\otimes_A R)$-submodule  $\mathcal{E} \subseteq \DiffR^n(M,N)$ such that the equality
		$$
		\Sol\big(\mathcal{E}\big) = \Sol\left(\mathcal{E}^\prime\right) \cap M
		$$
		holds.
	\end{enumerate}
	\begin{proof}
		$(i)$ 
		From \autoref{rem_conseq_essent_finite_type}$(i)$ we have that $k(\pp)$ is a finitely generated field over $\Quot(A)$.
		If $\iota:\Quot(A) \hookrightarrow k(\pp)$ is a separable extension, then we can choose elements $v_1,\ldots,v_d \in k(\pp
		)$ that form a separating transcendence basis over $\Quot(A)$ (see, e.g., \cite[Theorem 26.2]{MATSUMURA}), i.e., $v_1,\ldots,v_d$ are algebraically independent over $\Quot(A)$ and  $\iota\left(\Quot(A)\right)(v_1,\ldots,v_d) \subseteq k(\pp)$ is a separable algebraic extension.
		
		Let $x_1, \ldots, x_d\in R_\pp$ be elements such that $v_i \in k(\pp)$ is the residue class of $x_i$.
		Consider the $\Quot(A)$-algebra homomorphism 
		\begin{align}
			\label{eq_map_choose_transc_basis}
			\begin{split}
			\Quot(A)[x_1,\ldots,x_d] &\longrightarrow k(\pp) = R_\pp/\pp R_\pp\\
			x_i &\mapsto v_i
			\end{split}			
		\end{align}
		obtained as a composition of the canonical maps $\Quot(A)[x_1,\ldots,x_d] \hookrightarrow R_\pp$ and $R_\pp \twoheadrightarrow R_\pp/\pp R_\pp$.
		Since $v_1,\ldots,v_d$ are algebraically independent over $\Quot(A)$, the above map \autoref{eq_map_choose_transc_basis} is injective, and so all the non-zero elements of $\Quot(A)[x_1,\ldots,x_d]$ do not belong to the maximal ideal $\pp R_\pp$ of the local ring $R_\pp$.
		Therefore, every non-zero element of $\Quot(A)[x_1,\ldots,x_d]$ is a unit in $R_\pp$, and so the result follows by taking $\KK = \Quot(A)(x_1,\ldots,x_d) \subseteq R_\pp$.
		
		$(ii)$ From \autoref{rem_conseq_essent_finite_type}$(iii)$ and \autoref{lem_gen_props_diff_ops}$(iv)$, we have a canonical map 
		$$
		\Psi_1 : \DiffR^n(M, N) \rightarrow  R_\pp \otimes_R \DiffR^n(M,N) \cong  \Diff_{R_\pp/\Quot(A)}^n\left(M_\pp, N_\pp\right)
		$$
		that corresponds with localizing at the multiplicatively closed subset $R \setminus \pp$. 
		By using \autoref{lem_gen_props_diff_ops}$(ii)$, there is a canonical inclusion 
		$$
		\Psi_2 : \Diff_{R_\pp/\KK}^n\left(M_\pp, N_\pp\right) \hookrightarrow \Diff_{R_\pp/\Quot(A)}^n\left(M_\pp, N_\pp\right).
		$$
				
		Let $\mathcal{E}^\prime \subseteq \Diff_{R_\pp/\KK}^n\left(M_\pp,N_\pp\right)$ be an $\left(R_\pp\otimes_\KK R_\pp\right)$-submodule. 
		Then, we set 
		$$
		\mathcal{E} =    \Psi_1^{-1} \left( \Psi_2\left(\mathcal{E}^\prime\right)\right).
		$$
		It is clear that $\mathcal{E}$ is an $(R \otimes_A R)$-submodule of $\DiffR^n(M,N)$.
				
		Since $\Psi_1$ is just the localization map $\DiffR^n(M, N) \rightarrow  R_\pp \otimes_R \DiffR^n(M,N)$, for any $\delta^\prime \in \mathcal{E}^\prime$ there exists an element $r \in R \setminus \pp$ such that $\Psi_1^{-1} \left( \Psi_2 \left(\{r\delta^\prime\}\right)\right) \neq \emptyset$.
		For all $\delta^\prime \in \mathcal{E}^\prime$ and $\delta \in \Psi_1^{-1} \left( \Psi_2 \left( \{\delta^\prime\}\right)\right)$ we have the following commutative diagram (in principle, of $A$-homomorphisms)
		\begin{center}
		\begin{tikzpicture}[baseline=(current  bounding  box.center)]
		\matrix (m) [matrix of math nodes,row sep=2em,column sep=10em,minimum width=2em, text height=1.5ex, text depth=0.25ex]
		{
			M & N \\
			M_\pp &  N_\pp \\
		};
		\path[-stealth]
		(m-1-1) edge node [above]	{$\delta$} (m-1-2)
		(m-1-1) edge (m-2-1)
		(m-1-2) edge (m-2-2)
		(m-2-1) edge node [above]	{$\delta^\prime$} (m-2-2)
		;				
		\end{tikzpicture}	
		\end{center}
		and from the fact that $N$ is $R/\pp$ torsion-free, it follows that  $N \rightarrow N_\pp$ is injective and so $\Ker(\delta) = \Ker(\delta^\prime) \cap M$. 
		Hence, it is clear that $\Sol(\mathcal{E}) \supseteq \Sol(\mathcal{E}^\prime
		) \cap M$.
		Conversely, for any $\delta^\prime \in \mathcal{E}^\prime$, we can choose $r \in R \setminus \pp$ and $\delta \in \Psi_1^{-1} \left( \Psi_2 \left(\{r\delta^\prime\}\right)\right)$, then we  obtain $\Ker(\delta) = \Ker(r\delta^\prime) \cap M  = \Ker(\delta^\prime) \cap M$; so, $\Sol(\mathcal{E}) \subseteq \Sol(\mathcal{E}^\prime) \cap M$.
		Therefore, we have that $\Sol\big(\mathcal{E}\big) = \Sol\left(\mathcal{E}^\prime\right) \cap M$.
	\end{proof}
\end{proposition}

The lemma below deals with the injectivity of certain maps from modules over an Artinian local algebra, and it will be an important basic tool.

\begin{lemma}
	\label{lem_injective_map_tensor_prods}
	Let $\KK$ be a field, $S$ be an Artinian local $\KK$-algebra and $T$ be a $\KK$-algebra. 
	Let $M$ be an $S$-module.
	Then, for any $\mathcal{Q} \in \Spec\left(T \otimes_\KK S\right)$ prime ideal in $T \otimes_\KK S$ the canonical map
	$$
	M \;\rightarrow\; {\left(T \otimes_\KK M\right)}_\mathcal{Q}, \qquad
	m \in M \;\mapsto\; \frac{1 \otimes_\KK m}{1} \in {\left(T \otimes_\KK M\right)}_\mathcal{Q}
	$$
	is injective.
	\begin{proof}
		Fix any $\mathcal{Q} \in \Spec\left(T \otimes_\KK S\right)$ and suppose that $\mm$ is the unique prime ideal in $S$.
		Since $(S,\mm)$ is an Artinian local ring, we have that ${\left(1 \otimes_\KK \mm\right)}^k = 1 \otimes_\KK \mm^k = 0$ for some $k > 0$, and this implies that $1 \otimes_\KK \mm \subseteq \mathcal{Q}$.		
		Hence, $S \rightarrow {\left(T \otimes_\KK S\right)}_\mathcal{Q}$  is a local ring homomorphism, and being flat, it is then faithfully flat.  
		
		Therefore, $M \rightarrow {\left(T \otimes_\KK S\right)}_\mathcal{Q} \otimes_S M \cong  {\left(T \otimes_\KK M\right)}_\mathcal{Q}$ is injective.
	\end{proof}	
\end{lemma}

Next we divide the section into two separate subsections. 
For simplicity and because of the main interest in the case of ideals, we first treat primary ideals and then we concentrate on primary submodules.

\subsection{Primary ideals}
In the present subsection, we characterize primary ideals as solution sets of certain differential operators.

The following lemma contains some useful translations to be used later.

\begin{lemma}
	\label{lem_diff_ops_correspon_to_ideals}
	Let $\pp \in \Spec(R)$ be a prime ideal in $R$ and $\KK \subseteq R_\pp$ be a field.
	Assume that $k(\pp) \otimes_\KK R_\pp$ is a Noetherian ring.
	Then, the following statements hold:
	\begin{enumerate}[(i)]
		\item There is an isomorphism of $(R_\pp \otimes_\KK R_\pp)$-modules 
		$$
		\Diff_{R_\pp/\KK}^n(R_\pp,k(\pp)) \cong \Hom_{k(\pp)}\left(\frac{k(\pp) \otimes_\KK R_\pp}{\MM^{n+1}}, k(\pp)\right),
		$$
		where $\MM$ is the ideal given as the kernel of the canonical map 
		$$
		k(\pp) \otimes_\KK R_\pp \rightarrow k(\pp), \qquad k \otimes_\KK r \mapsto k\overline{r},
		$$ 
		and $\overline{r}$ denotes the residue class of $r$.		
		\item $\Diff_{R_\pp/\KK}^n(R_\pp,k(\pp))$ is a finite dimensional vector space over $k(\pp)$.
		\item There is a bijective correspondence between $\MM$-primary ideals $\mathcal{N} \subset k(\pp) \otimes_\KK R_\pp$ containing $\MM^{n+1}$ and $(R_\pp \otimes_\KK R_\pp)$-submodules of $\Diff_{R_\pp/\KK}^n\left(R_\pp,k(\pp)\right)$ induced by the isomorphism in part (i) and the mapping
		$$
		\mathcal{N} \supseteq \MM^{n+1} \;\mapsto\; \Hom_{k(\pp)}\left(\frac{k(\pp) \otimes_\KK R_\pp}{\mathcal{N}}, k(\pp)\right) \; \subseteq \; \Hom_{k(\pp)}\left(\frac{k(\pp) \otimes_\KK R_\pp}{\MM^{n+1}}, k(\pp)\right).
		$$
		\item Under the correspondence of part $(iii)$, if $\mathcal{E} \subseteq \Diff_{R_\pp/\KK}^n(R_\pp,k(\pp))$ is determined by an $\MM$-primary ideal $\mathcal{N} \supseteq \MM^{n+1}$, then we obtain the equality 
		$$
		\Sol(\mathcal{E}) = \mathcal{N} \cap R_\pp,
		$$
		where $\mathcal{N} \cap R_\pp$ denotes the contraction of 
		$\mathcal{N}$ under the canonical inclusion 
		$$
		R_\pp \cong 1 \otimes_\KK R_\pp \hookrightarrow k(\pp) \otimes_\KK R_\pp.
		$$
	\end{enumerate}
	\begin{proof}
		$(i)$ Using \autoref{prop_represen_diff_opp} and the Hom-tensor adjunction, we obtain the following isomorphisms 
		\begin{align*}		
		\Diff_{R_\pp/\KK}^n(R_\pp, k(\pp)) &\cong \Hom_{R_\pp}\left(P_{R_\pp/ \KK}^n, k(\pp)\right) \\
		&\cong \Hom_{R_\pp}\left(P_{R_\pp/ \KK}^n, \Hom_{k(\pp)}\left(k(\pp), k(\pp)\right)\right)\\
		&\cong \Hom_{k(\pp)}\left(k(\pp) \otimes_{R_\pp} P_{R_\pp/ \KK}^n, k(\pp)\right).	
		\end{align*}
		Finally, it is clear that 
		$$k(\pp) \otimes_{R_\pp} P_{R_\pp/\KK}^n \cong \frac{k(\pp) \otimes_\KK R_\pp}{\MM^{n+1}}.
		$$
		
		$(ii)$ Since $k(\pp) \otimes_\KK R_\pp$ is a Noetherian ring and $\MM \subset k(\pp) \otimes_\KK R_\pp$ is a maximal ideal, then $\frac{k(\pp) \otimes_\KK R_\pp}{\MM^{n+1}}$ is an Artinian local $k(\pp)$-algebra with residue field $k(\pp)$, and so it follows that $\frac{k(\pp) \otimes_\KK R_\pp}{\MM^{n+1}}$ is also a finite dimensional vector space over $k(\pp)$.
		Hence, the result is obtained from part $(i)$.
		
		$(iii)$
		By using part $(ii)$, the functor $\Hom_{k(\pp)}\left(\bullet, k(\pp)\right)$ gives a bijective correspondence between quotients of 
		$$
		\frac{k(\pp) \otimes_\KK R_\pp}{\MM^{n+1}} 
		$$
		by an ideal $\overline{\mathcal{N}} \subset \frac{k(\pp) \otimes_\KK R_\pp}{\MM^{n+1}} $ and $(R_\pp \otimes_\KK R_\pp)$-submodules of
		$$
		\Hom_{k(\pp)}\left(\frac{k(\pp) \otimes_\KK R_\pp}{\MM^{n+1}}, k(\pp)\right).
		$$
		Therefore, the statement follows from the isomorphism of part $(i)$.
		
		$(iv)$  From part $(iii)$, let $\mathcal{E} \subseteq \Diff_{R_\pp/\KK}^n\left(R_\pp, k(\pp)\right)$ such that 
		$$
		\mathcal{E} \cong \Hom_{k(\pp)}\left(\frac{k(\pp) \otimes_\KK R_\pp}{\mathcal{N}}, k(\pp)\right).
		$$
		Any $\delta \in \mathcal{E}$ as an element inside $\Diff_{R_\pp/\KK}^n(R_\pp,k(\pp)) \subseteq \Hom_\KK(R_\pp,k(\pp))$, via the isomorphism of part $(i)$, is given as $\delta(f) = \epsilon\left(\overline{1 \otimes_\KK f}\right)$ for $f \in R_\pp$, where  
		$$
		\epsilon \in \Hom_{k(\pp)}\left(\frac{k(\pp) \otimes_\KK R_\pp}{\mathcal{N}}, k(\pp)\right)
		$$ 
		(again, this follows from the Hom-tensor adjunction).
		So, for any $\delta \in \mathcal{E}$ we have that $\Ker(\delta) \supseteq \mathcal{N} \cap R_\pp$.
		In the other direction, for any $g \in R_\pp \setminus\ \left(\mathcal{N} \cap R_\pp\right)$, we can define a $k(\pp)$-linear map $\epsilon_g \in \Hom_{k(\pp)}\left(\frac{k(\pp) \otimes_\KK R_\pp}{\mathcal{N}}, k(\pp)\right)$ such that $\epsilon_g\left(\overline{1 \otimes_\KK g}\right)=1$, and so we obtain $\delta_g \in \mathcal{E} \subseteq  \Diff_{R_\pp/\KK}^n(R_\pp,k(\pp))$, given as $\delta_g(f)=\epsilon_g\left(\overline{1 \otimes_\KK f}\right)$ for $f \in R_\pp$, such that $\delta_g(g)=1 \neq 0$; so, $\Sol(\mathcal{E}) \subseteq \mathcal{N} \cap R_\pp$.
		Therefore, it follows that $\Sol(\mathcal{E}) = \mathcal{N} \cap R_\pp$.
	\end{proof}
\end{lemma}

The following proposition, in a separable setting, gives an isomorphism between the associated graded rings of $\pp R_\pp$ and $\MM$, where $\pp \in \Spec(R)$ and $\MM$ is the maximal ideal in \autoref{lem_diff_ops_correspon_to_ideals}.
The result of this proposition can also be found in \cite[Proposition 4.1]{BRUMFIEL_DIFF_PRIM}, but note that we provide a slightly different proof that depends upon \autoref{lem_injective_map_tensor_prods}.

\begin{proposition}[{\cite[Proposition 4.1]{BRUMFIEL_DIFF_PRIM}}]
	\label{prop_ass_gr_rings}
	Let $\pp \in \Spec(R)$ be a prime ideal in $R$, and $\KK \subseteq R_\pp$ be a field such that $\KK \hookrightarrow k(\pp)$ is a separable algebraic field extension.
	Let $\MM$ be the kernel of the canonical map $k(\pp) \otimes_\KK R_\pp \rightarrow k(\pp)$.
	Then, we have the following isomorphisms: 
	\begin{align*}
		&(i) \qquad
		\gr_{\pp R_\pp}\left(R_\pp\right) = \bigoplus_{n=0}^\infty \pp^nR_\pp/\pp^{n+1}R_\pp \;\;\xrightarrow{\cong}\;\;\gr_{\MM}\left(k(\pp) \otimes_\KK R_\pp\right) = \bigoplus_{n=0}^\infty \MM^n/\MM^{n+1}.\\
		&(ii) \qquad R_\pp/\pp^{n}R_\pp \;\;\xrightarrow{\cong}\;\; \left(k(\pp) \otimes_\KK R_\pp\right)/\MM^n \quad \text{ for all } n \ge 1.
	\end{align*}
 	\begin{proof}
 		$(i)$
		Note that $\MM$	is generated by elements of the form $1 \otimes_\KK r - \overline{r} \otimes_\KK 1$, where $r \in R_\pp$ and $\overline{r} \in k(\pp)$ is the residue class of $r$.
		Indeed, for any $\sum_{i=1}^k \alpha_i \otimes_\KK r_i$ where $\alpha_i \in k(\pp)$, $r_i \in R_\pp$ and $\sum_{i=1}^k\alpha_i\overline{r_i}=0$, we can write 
		$$
		\sum_{i=1}^k \alpha_i \otimes_\KK r_i = \sum_{i=1}^k \left(\alpha_i \otimes_\KK 1\right)\left(1 \otimes_\KK r_i - \overline{r_i} \otimes_\KK 1\right).
		$$
		
		We denote by $\left(1 \otimes_\KK \pp^n R_\pp\right) \subset k(\pp) \otimes_\KK R_\pp$ the ideal generated by the elements of $1 \otimes_\KK \pp^n R_\pp$, that is, the extension of $\pp^n R_\pp$ under the canonical map 
		$$
		R_\pp \cong 1 \otimes_\KK R_\pp \hookrightarrow k(\pp) \otimes_\KK R_\pp.
		$$
		Since $\left(1 \otimes_\KK \pp R_\pp\right) \subseteq \MM$, there are well-defined maps $\pp^nR_\pp/\pp^{n+1}R_\pp \rightarrow \MM^n/\MM^{n+1}$.
		Also, the case $n=0$ is clear because the map $k(\pp) \cong R_\pp/\pp R_\pp \rightarrow \left(k(\pp) \otimes_\KK R_\pp\right)/\MM \cong k(\pp)$ is an isomorphism.
		
		First, we prove that $\pp^nR_\pp/\pp^{n+1}R_\pp \rightarrow \MM^n/\MM^{n+1}$ is surjective. 
		But, actually it is enough to show that $\pp R_\pp/\pp^{2}R_\pp \rightarrow \MM/\MM^{2}$ is surjective.
		Let $r \in R_\pp$ and $f(x) \in \KK[x] \subseteq R_\pp[x]$ be the minimal polynomial of $\overline{r}$ over $\KK$.
		From $f(\overline{r})=0$ and $f^\prime(\overline{r})\neq 0$, we obtain that $f(r) \in \pp R_\pp$ and $f^\prime(r) \not\in \pp R_\pp$ is a unit in $R_\pp$.
		By taking the Taylor expansion
		$$
		0 = f(\overline{r} \otimes_\KK 1) = f(1\otimes_\KK r) + f^\prime(1 \otimes_\KK r)\left(\overline{r} \otimes_\KK 1 - 1\otimes_\KK r\right) \left(\text{mod}\; \MM^2\right), 
		$$
		it follows 
		$$
		1 \otimes_\KK r - \overline{r} \otimes_\KK 1 = \frac{f(1\otimes_\KK r)}{f^\prime(1\otimes_\KK r)} \left(\text{mod}\; \MM^2\right).
		$$
		So, the map $\pp R_\pp/\pp^{2}R_\pp \rightarrow \MM/\MM^{2}$ is surjective. 
		
		Now, we prove that $R_\pp/\pp^n R_\pp \rightarrow \left(k(\pp) \otimes_\KK R_\pp\right)/\MM^n$ is injective.
		From \autoref{lem_injective_map_tensor_prods} we have that the canonical map
		\begin{equation}
			\label{eq_injection_powers_pp}
			R_\pp/\pp^nR_\pp \rightarrow {\left(k(\pp) \otimes_\KK R_\pp/\pp^nR_\pp\right)}_\MM \cong {\left(\frac{k(\pp) \otimes_\KK R_\pp}{\left(1\otimes_\KK \pp^nR_\pp\right)}\right)}_\MM
		\end{equation}		
		is injective.
		The surjectivity of the map $\pp^{n} R_\pp \rightarrow \MM^n/\MM^{n+1}$ yields 
		$$
		\MM^{n}=\left(1\otimes_\KK \pp^nR_\pp\right)  + \MM \cdot \MM^n,
		$$ 
		then Nakayama's lemma applied in the local ring ${\left(k(\pp) \otimes_\KK R_\pp\right)}_\MM$ gives us 
		\begin{equation}
			\label{eq_isom_powers_pp}
			\left(1\otimes_\KK \pp^nR_\pp\right){\left(k(\pp) \otimes_\KK R_\pp\right)}_\MM = \MM^n{\left(k(\pp) \otimes_\KK R_\pp\right)}_\MM.
		\end{equation}		
		By summing up \autoref{eq_injection_powers_pp} and \autoref{eq_isom_powers_pp} we obtain that
		$$ 
		R_\pp/\pp^nR_\pp \rightarrow {\left(\left(k(\pp) \otimes_\KK R_\pp\right)/\MM^n\right)}_\MM = \left(k(\pp) \otimes_\KK R_\pp\right)/\MM^n
		$$
		is injective, as required. 
		Therefore, the first isomorphism follows. 
		
		$(ii)$ It is obtained inductively from the exact sequences 
		\begin{align*}
			0 \rightarrow \pp^n&R_\pp/\pp^{n+1}R_\pp \rightarrow R_\pp/\pp^{n+1}R_\pp \rightarrow R_\pp/\pp^{n}R_\pp \rightarrow 0 \\
			0 \rightarrow \MM^n/\MM^{n+1}& \rightarrow \left(k(\pp) \otimes_\KK R_\pp\right)/\MM^{n+1} \rightarrow \left(k(\pp) \otimes_\KK R_\pp\right)/\MM^{n} \rightarrow 0
		\end{align*}
		and part $(i)$.
 	\end{proof}
\end{proposition}

\begin{remark}
		\label{rem:finite_gen_Noeth_ops}
		(i) In the setting of  \autoref{lem_diff_ops_correspon_to_ideals} we have that $\Diff_{R_\pp/\KK}^n(R_\pp,k(\pp))$ is a finitely generated $R_\pp$-module.
		If $\{\delta_1, \ldots,\delta_q\}$ is a set of generators of $\mathcal{E} \subseteq \Diff_{R_\pp/\KK}^n(R_\pp,k(\pp))$ as an $R_\pp$-module, then $\Sol(\mathcal{E}) = \Sol\left(\{\delta_1, \ldots,\delta_q\}\right)$.
		Indeed, for any $\delta \in \mathcal{E}$, we can then write $\delta = r_1\delta_1 + \cdots + r_q\delta_q$ for some $r_i \in R_\pp$;
		so $\delta_i(r)=0$ implies $\delta(r)=0$ for every $r \in R_\pp$.
		
		\smallskip
		\noindent
		(ii) If $R$ is essentially of finite type over $A$ and $N$ is a finitely generated $R$-module,  \autoref{rem_conseq_essent_finite_type} and \autoref{prop_represen_diff_opp} yield that $\DiffR^n(M,N)$ is a finitely generated $R$-module.
		If $\{\delta_1, \ldots,\delta_q\}$ is a set of generators of $\mathcal{E} \subseteq \DiffR^n(M,N)$ as an $R$-module, then $\Sol(\mathcal{E}) = \Sol\left(\{\delta_1, \ldots,\delta_q\}\right)$.
\end{remark}

We are now ready for the main result about primary ideals. 
The following theorem and \autoref{cor_prim_ideals_sol_Diff_R_R} below generalize the main results of \cite{BRUMFIEL_DIFF_PRIM} and \cite{OBERST_NOETH_OPS} regarding primary ideals.

\begin{theorem}
	\label{thm_Pimary_ideals_as_Sol}
	Assume \autoref{setup_primary_submod}. 
	Let $\pp \in \Spec(R)$ be a prime ideal in $R$ such that $\pp \cap A = 0$.
	Then, the following statements hold:
	\begin{enumerate}[(i)]
		\item Suppose that $k(\pp) \otimes_{\Quot(A)} R_\pp$ is a Noetherian ring.
		If $I \subset R$ is a $\pp$-primary ideal in $R$, then there exists an $(R \otimes_A R)$-submodule $\mathcal{E} \subseteq \DiffR(R,k(\pp))$ such that 
		$$
		I = \Sol(\mathcal{E})
		$$
		and $\mathcal{E}$ is finitely generated as an $R_\pp$-module.
		\item Suppose that $R$ is essentially of finite type over $A$ and $N$ is a finitely generated torsion-free module over $R/\pp$.
		\begin{enumerate}[(a)]
			\item If $I \subset R$ is a $\pp$-primary ideal in $R$, then there exists an $(R \otimes_A R)$-submodule $\mathcal{E} \subseteq \DiffR(R,N)$ such that 
			$$
			I = \Sol(\mathcal{E})
			$$
			and $\mathcal{E}$ is finitely generated as an $R$-module.
			\item If $\Quot(A) \hookrightarrow k(\pp)$ is a separable field extension, which holds whenever $\Quot(A)$ is perfect, then for any $\pp$-primary ideal $I \subset R$ containing $\pp^{n+1}$, there exists an $(R \otimes_A R)$-submodule $\mathcal{E} \subseteq \DiffR^n(R,N)$ such that 
			$$
			I = \Sol(\mathcal{E}).
			$$
		\end{enumerate}				
	\end{enumerate}
	\begin{proof}
		For notational purposes, set $Q = \Quot(A)$.
		
		$(i)$ Since $k(\pp)$ is clearly an $R_\pp$-module, \autoref{lem_gen_props_diff_ops}$(iii)$ gives the following isomorphism 
		$$
		\Psi : \DiffR^n(R, k(\pp)) \xrightarrow{\cong} \Diff_{R_\pp/Q}^n(R_\pp, k(\pp)),
		$$
		and for any $\delta \in \DiffR^n(R, k(\pp))$ we have that $\Ker(\delta) = \Ker\left(\Psi(\delta)\right) \cap R$ (once again, this follows from the Hom-tensor adjunction; see \autoref{eq_diff_op_N_is_localized}).
		Therefore, using the fact that $I$ is a $\pp$-primary ideal, it suffices to find an $(R_\pp \otimes_Q R_\pp)$-submodule $\mathcal{E}^\prime \subseteq \Diff_{R_\pp/Q}^n\left(R_\pp, k(\pp)\right)$ such that $IR_\pp=\Sol\left(\mathcal{E}^\prime\right)$.
		Then, we can take $\mathcal{E} = \Psi^{-1}\left(\mathcal{E}^\prime\right)$ and the claim about $\mathcal{E}$ being finitely generated as an $R_\pp$-module is clear from \autoref{rem:finite_gen_Noeth_ops}.
		
		By using the same argument of \autoref{prop_reduction_to_zero_dim}$(i)$ we choose a field $Q \subseteq \KK \subseteq R_\pp$ such that $\KK \hookrightarrow k(\pp)$ is an algebraic extension (again, we can take a transcendence basis of $k(\pp)$ over $Q$ and then pull it back to elements of $R_\pp$).
		Then, under the assumption that $\KK \hookrightarrow k(\pp)$ is an algebraic extension, one has that 
		\begin{equation*}
			R_\pp/IR_\pp \hookrightarrow k(\pp) \otimes_\KK R_\pp/IR_\pp
		\end{equation*}
		is an integral extension and so $\dim\left(k(\pp) \otimes_\KK R_\pp/IR_\pp\right)=\dim\left(R_\pp/IR_\pp\right)=0$ (see, e.g., \cite[Lemma 2.4]{SHARP_TENSOR_PROD_FIELDS}).
		Since  $k(\pp) \otimes_{Q} R_\pp$ is Noetherian, it follows that $k(\pp) \otimes_\KK R_\pp/IR_\pp$ is an Artinian ring.
		The canonical inclusion 
		$$
		\Diff_{R_\pp/\KK}^n(R_\pp, k(\pp)) \hookrightarrow \Diff_{R_\pp/Q}^n(R_\pp, k(\pp))
		$$
		(see \autoref{lem_gen_props_diff_ops}$(ii)$) yields that it is enough to find an $(R_\pp \otimes_\KK R_\pp)$-submodule $\mathcal{E}^\prime \subseteq \Diff_{R_\pp/\KK}^n\left(R_\pp, k(\pp)\right)$ such that $IR_\pp=\Sol\left(\mathcal{E}^\prime\right)$.
			
		Let $\MM \subset k(\pp) \otimes_\KK R_\pp$ be the ideal given as the kernel of the canonical map $k(\pp) \otimes_\KK R_\pp \rightarrow k(\pp)$.
		By using \autoref{lem_diff_ops_correspon_to_ideals}, if  we show the existence of an $\MM$-primary ideal $\mathcal{N} \subset k(\pp) \otimes_\KK R_\pp$ such that $IR_\pp = \mathcal{N} \cap R_\pp$, then the result would be obtained.
		Now, we proceed to find such an ideal $\mathcal{N}$.
		
		From \autoref{lem_injective_map_tensor_prods} the canonical map
		$$
		\alpha : R_\pp/IR_\pp \rightarrow {\left(k(\pp) \otimes_\KK R_\pp/IR_\pp\right)}_{\MM}
		$$
		is injective.
		Denote by $\left(1 \otimes_\KK IR_\pp\right) \subset k(\pp) \otimes_\KK R_\pp$ the ideal generated by the elements of $1 \otimes_\KK IR_\pp$.
		We set $\mathcal{N} \subset k(\pp) \otimes_\KK R_\pp$ to be the kernel of the map 
		$$
		k(\pp) \otimes_\KK R_\pp \rightarrow \frac{k(\pp) \otimes_\KK R_\pp}{\left(1 \otimes_\KK IR_\pp\right)} \cong  {k(\pp) \otimes_\KK R_\pp/IR_\pp} \rightarrow {\left(k(\pp) \otimes_\KK R_\pp/IR_\pp\right)}_{\MM}.
		$$
		Since ${\left(k(\pp) \otimes_\KK R_\pp/IR_\pp\right)}_{\MM}$ is an Artinian local ring, $\mathcal{N}$ is an $\MM$-primary ideal. 
		From the following commutative diagram
		\begin{center}
		\begin{tikzpicture}[baseline=(current  bounding  box.center)]
		\matrix (m) [matrix of math nodes,row sep=2em,column sep=4em,minimum width=2em, text height=1.5ex, text depth=0.25ex]
		{
			k(\pp) \otimes_\KK R_\pp & {k(\pp) \otimes_\KK R_\pp/IR_\pp} & {\left(k(\pp) \otimes_\KK R_\pp/IR_\pp\right)}_{\MM} \\
			R_\pp &  R_\pp/IR_\pp \\
		};
		\path[-stealth]		
		(m-1-1) edge (m-1-2)
		(m-1-2) edge (m-1-3)
		(m-2-2) edge (m-1-2)
		(m-2-1) edge (m-1-1)
		(m-2-1) edge (m-2-2)
		(m-2-2) edge node [below] {$\alpha$} (m-1-3)
		;				
		\end{tikzpicture}	
		\end{center}		 
		and the injectivity of $\alpha$, we obtain that $IR_\pp = \mathcal{N} \cap R_\pp$.
		So, the proof of this part follows.
		
		$(ii.a)$ From \autoref{prop_reduction_to_zero_dim}$(ii)$ we have that for any   $\left(R_\pp\otimes_Q R_\pp\right)$-submodule $\mathcal{E}^\prime \subseteq \Diff_{R_\pp/Q}^n\left(R_\pp, N_\pp\right)$, there exists an  $(R\otimes_A R)$-submodule  $\mathcal{E} \subseteq \DiffR^n(R,N)$ that satisfies the equality
		$
		\Sol\big(\mathcal{E}\big) = \Sol\left(\mathcal{E}^\prime\right) \cap R.
		$
		
		Hence, since $I$ is a $\pp$-primary ideal, it is enough to find an  $(R_\pp \otimes_Q R_\pp)$-submodule $\mathcal{E}^\prime \subseteq \Diff_{R_\pp/Q}^n\left(R_\pp, N_\pp\right)$ such that $IR_\pp=\Sol\left(\mathcal{E}^\prime\right)$.
		Since $N$ is a finitely generated $R/\pp$-module, note that $N_\pp$ is a finite dimensional vector space over $k(\pp)$, say $r = \dim_{k(\pp)}\left(N_\pp\right)$, thus \autoref{lem_gen_props_diff_ops}$(vi)$ yields the isomorphism 
		$$
		\Diff_{R_\pp/Q}^n\left(R_\pp, N_\pp\right) \cong {\left(\Diff_{R_\pp/Q}^n\left(R_\pp, k(\pp)\right)\right)}^r.
		$$
		Therefore, it is enough to consider the case where $N_\pp = k(\pp)$, and so the result follows from part $(i)$.
		
		$(ii.b)$ First, if $Q$ is perfect, then $Q \hookrightarrow k(\pp)$ is a separable extension (see, e.g., \cite[Theorem 26.3]{MATSUMURA}).
		From \autoref{prop_reduction_to_zero_dim}$(i)$, since $Q \hookrightarrow k(\pp)$ is assumed to be separable, there exists a field $Q \subseteq \KK \subseteq R_\pp$ such that $\KK \hookrightarrow k(\pp)$ is a separable finite extension. 
		By following the same steps of part $(ii.a)$, from \autoref{prop_reduction_to_zero_dim}$(ii)$, \autoref{lem_gen_props_diff_ops}$(vi)$ and \autoref{lem_diff_ops_correspon_to_ideals}, now it is enough to find an $\MM$-primary ideal $\mathcal{N} \subset k(\pp) \otimes_\KK R_\pp$ such that $IR_\pp = \mathcal{N} \cap R_\pp$ and $\mathcal{N} \supseteq \MM^{n+1}$.
		
		Under the assumption that $\KK \hookrightarrow k(\pp)$ is a separable finite extension, \autoref{prop_ass_gr_rings}$(ii)$ yields the canonical isomorphism 
		$$
		\Phi : R_\pp/\pp^{n+1}R_\pp \;\xrightarrow{\cong}\; \left(k(\pp) \otimes_\KK R_\pp\right)/\MM^{n+1}.
		$$
		Therefore, since $I \supseteq \pp^{n+1}$, then the result follows by taking $\mathcal{N}$ as the preimage in $k(\pp) \otimes_\KK R_\pp$ of the ideal $\Phi\left(IR_\pp/\pp^{n+1}R_\pp\right) \subset \left(k(\pp) \otimes_\KK R_\pp\right)/\MM^{n+1}$.
	\end{proof}
\end{theorem}

A natural question after the previous theorem is whether the differential operators inside $\DiffR(R, R)$ are enough to characterize primary ideals; we shall see that under certain smooth settings it is actually possible. 
Before, we note the following result that will allow us to lift differential operators.

\begin{proposition}
	\label{prop_exactness_diff_op}
	Assume \autoref{setup_primary_submod} with the additional condition that $R$ is formally smooth and essentially of finite type over $A$.
	Let $F$ be a finitely generated free $R$-module.
	Then, $\DiffR^n(F, \bullet)$ is a covariant exact functor.
	\begin{proof}
		From \autoref{prop_represen_diff_opp} it is enough to show that $\Princ^n(F)\cong \Princ^n \otimes_R F$ is a projective $R$-module; equivalently, we can show that $\Princ^n$ is a projective $R$-module. 
		But, by using \cite[Proposition 16.10.2]{EGAIV_IV} and \cite[D\'efinition 16.10.1]{EGAIV_IV} we obtain that $\Princ^n$ is a projective $R$-module.
	\end{proof}
\end{proposition}

Under the assumption of $R$ being formally smooth and essentially of finite type over $A$, we show that the same results of \autoref{thm_Pimary_ideals_as_Sol}$(ii)$ also hold by using differential operators in $\DiffR(R,R)$.

\begin{corollary}
	\label{cor_prim_ideals_sol_Diff_R_R}
	Assume \autoref{setup_primary_submod} with the additional condition that $R$ is formally smooth and essentially of finite type over $A$.
	Let $\pp \in \Spec(R)$ be a prime ideal in $R$ such that $\pp \cap A = 0$.
	Then, the following statements hold:
	\begin{enumerate}[(i)]
		\item If $I \subset R$ is a $\pp$-primary ideal in $R$, then there exists an $(R \otimes_A R)$-submodule $\mathcal{E} \subseteq \DiffR(R,R)$ such that 
		$$
		I = \big\lbrace f \in R \mid \delta(f) \in \pp \text{ for all } \delta \in \mathcal{E} \big\rbrace
		$$
		and $\mathcal{E}$ is finitely generated as an $R$-module.
		\item If $\Quot(A) \hookrightarrow k(\pp)$ is a separable field extension, which holds whenever $\Quot(A)$ is perfect, then for any $\pp$-primary ideal $I \subset R$ containing $\pp^{n+1}$, there exists an  $(R \otimes_A R)$-submodule $\mathcal{E} \subseteq \DiffR^n(R,R)$ such that 
		$$
		I = \big\lbrace f \in R \mid \delta(f) \in \pp \text{ for all } \delta \in \mathcal{E} \big\rbrace.
		$$
	\end{enumerate}
	\begin{proof}
		$(i)$ From \autoref{thm_Pimary_ideals_as_Sol}$(ii.a)$, let $\mathcal{E}^\prime \subseteq \DiffR^n(R,R/\pp)$ be an  $(R \otimes_A R)$-submodule such that 
		$
		I = \Sol(\mathcal{E}^\prime).
		$
		Since $R$ is formally smooth and essentially of finite type over $A$, \autoref{prop_exactness_diff_op} yields a canonical surjection 
		$$
		\Psi: \DiffR^n(R, R) \twoheadrightarrow \DiffR^n(R,R/\pp).
		$$
		So, the result follows by taking $\mathcal{E} = \Psi^{-1}\left(\mathcal{E}^\prime\right)$.
		
		$(ii)$ Follows identically to part $(i)$.
	\end{proof}
\end{corollary}

\subsection{Primary submodules}
Here we extend the results of the previous subsection to describe primary submodules.
The proofs in this subsection will be relatively easy adaptations.

The lemma below is a simple extension of \autoref{lem_diff_ops_correspon_to_ideals}.

\begin{lemma}
	\label{lem_diff_ops_correspon_to_submodules}
	Let $\pp \in \Spec(R)$ be a prime ideal in $R$ and $\KK \subseteq R_\pp$ be a field.
	Assume that $k(\pp) \otimes_\KK R_\pp$ is a Noetherian ring.
	Let $M$ be a finitely generated $R$-module.
	Then, the following statements hold:
	\begin{enumerate}[(i)]
		\item There is an isomorphism of $(R_\pp \otimes_\KK R_\pp)$-modules 
		$$
		\Diff_{R_\pp/\KK}^n(M_\pp,k(\pp)) \cong \Hom_{k(\pp)}\left(\frac{k(\pp) \otimes_\KK M_\pp}{\MM^{n+1}\left(k(\pp) \otimes_\KK M_\pp\right)}, k(\pp)\right),
		$$
		where, as in \autoref{lem_diff_ops_correspon_to_ideals}, $\MM$ is the kernel of the canonical map 
		$
		k(\pp) \otimes_\KK R_\pp \rightarrow k(\pp).
		$ 
		\item $\Diff_{R_\pp/\KK}^n(M_\pp,k(\pp))$ is a finite dimensional vector space over $k(\pp)$.
		\item There is a bijective correspondence between $\MM$-primary submodules 
		$
		\mathbb{V} \subset k(\pp) \otimes_\KK M_\pp
		$ 
		containing $\MM^{n+1}\left(k(\pp) \otimes_\KK M_\pp\right)$ and  $(R_\pp \otimes_\KK R_\pp)$-submodules of 
		$
		\Diff_{R_\pp/\KK}^n(M_\pp,k(\pp))
		$
		induced by the isomorphism in part (i) and the mapping
		$$
		\mathbb{V} \supseteq \MM^{n+1}\left(k(\pp) \otimes_\KK M_\pp\right) \;\mapsto\; \Hom_{k(\pp)}\left(\frac{k(\pp) \otimes_\KK M_\pp}{\mathbb{V}}, k(\pp)\right).
		$$
		\item Under the correspondence of part $(iii)$, if $\mathcal{E} \subseteq \Diff_{R_\pp/\KK}^n(M_\pp,k(\pp))$ is determined by an $\MM$-primary submodule $\mathbb{V} \supseteq \MM^{n+1}\left(k(\pp) \otimes_\KK M_\pp\right)$ then we obtain the equality 
		$$
		\Sol(\mathcal{E}) = \mathbb{V} \cap M_\pp,
		$$
		where $\mathbb{V} \cap M_\pp$ denotes the contraction of 
		$\mathbb{V}$ under the canonical inclusion 
		$$
		M_\pp \cong 1 \otimes_\KK M_\pp \hookrightarrow k(\pp) \otimes_\KK M_\pp.
		$$
	\end{enumerate}
	\begin{proof}
		$(i)$ We can use the same proof of \autoref{lem_diff_ops_correspon_to_ideals}$(i)$ by only noting the isomorphisms 
		\begin{align*}
			k(\pp) \otimes_{R_\pp} P_{R_\pp/\KK}^n(M_\pp) &\cong k(\pp) \otimes_{R_\pp} P_{R_\pp/\KK}^n \otimes_{R_\pp} M_\pp   \\
			&\cong \frac{k(\pp) \otimes_\KK R_\pp}{\MM^{n+1}} \otimes_{R_\pp} M_\pp\\
			&\cong \frac{k(\pp) \otimes_\KK M_\pp}{\MM^{n+1}\left(k(\pp) \otimes_\KK M_\pp\right)}
		\end{align*}
		(see \autoref{rem_isom_mod_princ_parts}).
		
		$(ii)$ Since $M$ is a finitely generated $R$-module,  
		$\frac{k(\pp) \otimes_\KK M_\pp}{\MM^{n+1}\left(k(\pp) \otimes_\KK M_\pp\right)}$ is a finite dimensional vector space over $k(\pp)$.

		$(iii)$ As in \autoref{lem_diff_ops_correspon_to_ideals}$(iii)$, the correspondence is given by the functor $\Hom_{k(\pp)}\left(\bullet, k(\pp)\right)$.
		
		$(iv)$  It follows verbatim to \autoref{lem_diff_ops_correspon_to_ideals}$(iv)$.
	\end{proof}
\end{lemma}

The following theorem is an extension of \autoref{thm_Pimary_ideals_as_Sol} to the case of primary submodules.
This theorem and \autoref{cor_prim_submods_sol_Diff_R_R} below generalize the main results of \cite{BRUMFIEL_DIFF_PRIM} and \cite{OBERST_NOETH_OPS} regarding primary submodules.

\begin{theorem}
	\label{thm_Pimary_submod_as_Sol}
	Assume \autoref{setup_primary_submod}. 
	Let $\pp \in \Spec(R)$ be a prime ideal in $R$ such that $\pp \cap A = 0$.
	Let $M$ be a finitely generated $R$-module.
	Then, the following statements hold:
	\begin{enumerate}[(i)]
		\item Suppose that $k(\pp) \otimes_{\Quot(A)} R_\pp$ is a Noetherian ring.
		If $U \subset M$ is a $\pp$-primary $R$-submodule, then there exists an  $(R \otimes_A R)$-submodule $\mathcal{E} \subseteq \DiffR(M,k(\pp))$ such that 
		$$
		U = \Sol(\mathcal{E})
		$$
		and $\mathcal{E}$ is finitely generated as an $R_\pp$-module.
		\item Suppose that $R$ is essentially of finite type over $A$ and $N$ is a finitely generated torsion-free module over $R/\pp$.
		\begin{enumerate}[(a)]
			\item If $U \subset M$ is a $\pp$-primary $R$-submodule, then there exists an  $(R \otimes_A R)$-submodule $\mathcal{E} \subseteq \DiffR(M,N)$ such that 
			$$
			U = \Sol(\mathcal{E})
			$$
			and $\mathcal{E}$ is finitely generated as an $R$-module.
			\item If $\Quot(A) \hookrightarrow k(\pp)$ is a separable field extension, which holds whenever $\Quot(A)$ is perfect, then for any $\pp$-primary $R$-submodule $U \subset M$ containing $\pp^{n+1}M$, there exists an  $(R \otimes_A R)$-submodule $\mathcal{E} \subseteq \DiffR^n(M,N)$ such that 
			$$
			U = \Sol(\mathcal{E}).
			$$
		\end{enumerate}			
	\end{enumerate}
	\begin{proof}
		The proof is completely similar to the one of \autoref{thm_Pimary_ideals_as_Sol}. Set $Q = \Quot(A)$.
		
		$(i)$ From \autoref{lem_gen_props_diff_ops}$(iii)$ we obtain the following isomorphism 
		$$
		\Psi : \DiffR^n(M, k(\pp)) \xrightarrow{\cong} \Diff_{R_\pp/Q}^n(M_\pp, k(\pp)).
		$$
		The claim about $\mathcal{E}$ being finitely generated as an $R_\pp$-module is also clear in this case (see \autoref{rem:finite_gen_Noeth_ops} and \autoref{lem_diff_ops_correspon_to_submodules}).
		
		As in \autoref{thm_Pimary_ideals_as_Sol}$(i)$, we can find a field $Q \subseteq \KK \subseteq R_\pp$ such that $\KK \hookrightarrow k(\pp)$ is an algebraic extension.
		Set $\mathfrak{b}=\Ann_{R_\pp}\left(M_\pp/U_\pp\right) \subset R_\pp$, since $U$ is a $\pp$-primary submodule, then $R_\pp/\mathfrak{b}$ is an Artinian ring, and, again, the integral extension
		$$
		R_\pp/\mathfrak{b} \hookrightarrow k(\pp) \otimes_\KK R_\pp/\mathfrak{b}
		$$
		yields that $k(\pp) \otimes_\KK R_\pp/\mathfrak{b}$ is an Artinian ring.		 
		Hence, $k(\pp) \otimes_\KK M_\pp/U_\pp$ is a module of finite length because it is finitely generated over $k(\pp) \otimes_\KK R_\pp/\mathfrak{b}$.
		
		Then, by following the same steps of \autoref{thm_Pimary_ideals_as_Sol}$(i)$, and now using \autoref{lem_diff_ops_correspon_to_submodules} instead of \autoref{lem_diff_ops_correspon_to_ideals}, it is enough to show the existence of an $\MM$-primary submodule $\mathbb{V} \subset k(\pp) \otimes_\KK M_\pp$ such that $U_\pp = \mathbb{V} \cap M_\pp$.

		From \autoref{lem_injective_map_tensor_prods} the canonical map
		$$
		\alpha: M_\pp/U_\pp \rightarrow {\left(k(\pp) \otimes_\KK M_\pp/U_\pp\right)}_{\MM}
		$$
		is injective.
		Denote by $\left(1 \otimes_\KK U_\pp\right) \subset k(\pp) \otimes_\KK M_\pp$ the submodule generated by the elements of $1 \otimes_\KK U_\pp$.
		We set $\mathbb{V}$ to be the kernel of the map 
		$$
		k(\pp) \otimes_\KK M_\pp \rightarrow \frac{k(\pp) \otimes_\KK M_\pp}{\left(1 \otimes_\KK U_\pp\right)} \cong  {k(\pp) \otimes_\KK M_\pp/U_\pp} \rightarrow {\left(k(\pp) \otimes_\KK M_\pp/U_\pp\right)}_{\MM}.
		$$
		Since ${\left(k(\pp) \otimes_\KK M_\pp/U_\pp\right)}_{\MM}$ is a module of finite length, it follows that $\mathbb{V}$ is an $\MM$-primary submodule.
		From the following commutative diagram
		\begin{center}
			\begin{tikzpicture}[baseline=(current  bounding  box.center)]
			\matrix (m) [matrix of math nodes,row sep=2em,column sep=4em,minimum width=2em, text height=1.5ex, text depth=0.25ex]
			{
				k(\pp) \otimes_\KK M_\pp & {k(\pp) \otimes_\KK M_\pp/U_\pp} & {\left(k(\pp) \otimes_\KK M_\pp/U_\pp\right)}_{\MM} \\
				M_\pp &  M_\pp/U_\pp \\
			};
			\path[-stealth]		
			(m-1-1) edge (m-1-2)
			(m-1-2) edge (m-1-3)
			(m-2-2) edge (m-1-2)
			(m-2-1) edge (m-1-1)
			(m-2-1) edge (m-2-2)
			(m-2-2) edge node [below] {$\alpha$} (m-1-3)
			;				
			\end{tikzpicture}	
		\end{center}		 
		and the injectivity of $\alpha$, we obtain that $U_\pp = \mathbb{V} \cap M_\pp$.
		So, the proof of this part also follows.
		
		$(ii.a)$ It follows verbatim to \autoref{thm_Pimary_ideals_as_Sol}$(ii.a)$, but now using the above part $(i)$ instead of \autoref{thm_Pimary_ideals_as_Sol}$(i)$.
		
		$(ii.b)$ As in \autoref{thm_Pimary_ideals_as_Sol}$(ii.b)$, we choose a field $Q \subseteq \KK \subseteq R_\pp$ such that $\KK \hookrightarrow k(\pp)$ is a separable finite extension, and here it is enough to find an $\MM$-primary submodule $\mathbb{V} \subset k(\pp) \otimes_\KK M_\pp$ such that $U_\pp = \mathbb{V} \cap M_\pp$ and $\mathbb{V} \supseteq \MM^{n+1}\left(k(\pp) \otimes_\KK M_\pp\right)$.
		
		By using \autoref{prop_ass_gr_rings}$(ii)$ and taking tensor product with $M_\pp$, we obtain the canonical isomorphism 
		\begin{align*}
			\Phi : M_\pp/\pp^{n+1}M_\pp \cong R_\pp/\pp^{n+1}R_\pp \otimes_{R_\pp} M_\pp \xrightarrow{\cong} &\left(k(\pp) \otimes_\KK R_\pp\right)/\MM^{n+1} \otimes_{R_\pp} M_\pp \\
			&\cong \frac{k(\pp) \otimes_\KK M_\pp}{\MM^{n+1}\left(k(\pp) \otimes_\KK M_\pp\right)}.
		\end{align*}
		Therefore, since $U \supseteq \pp^{n+1}M$, then the result follows by taking $\mathbb{V}$ as the preimage in $k(\pp) \otimes_\KK M_\pp$ of the submodule $\Phi\left(U_\pp/\pp^{n+1}M_\pp\right) \subset \left(k(\pp) \otimes_\KK M_\pp\right)/\MM^{n+1}\left(k(\pp) \otimes_\KK M_\pp\right)$.
	\end{proof}
\end{theorem}

Finally, we provide an extension of \autoref{cor_prim_ideals_sol_Diff_R_R} for primary submodules of a finitely generated free $R$-module.

\begin{corollary}
	\label{cor_prim_submods_sol_Diff_R_R}
	Assume \autoref{setup_primary_submod} with the additional condition that $R$ is formally smooth and essentially of finite type over $A$.
	Let $F$ be a finitely generated free $R$-module.
	Let $\pp \in \Spec(R)$ be a prime ideal in $R$ such that $\pp \cap A = 0$.
	Then, the following statements hold:
	\begin{enumerate}[(i)]
		\item If $U \subset F$ is a $\pp$-primary $R$-submodule, then there exists an  $(R \otimes_A R)$-submodule $\mathcal{E} \subseteq \DiffR(F,R)$ such that 
		$$
		U = \big\lbrace f \in F \mid \delta(f) \in \pp \text{ for all } \delta \in \mathcal{E} \big\rbrace
		$$
		and $\mathcal{E}$ is finitely generated as an $R$-module.
		\item If $\Quot(A) \hookrightarrow k(\pp)$ is a separable field extension, which holds whenever $\Quot(A)$ is perfect, then for any $\pp$-primary $R$-submodule $U \subset F$ containing $\pp^{n+1}F$, there exists an  $(R \otimes_A R)$-submodule $\mathcal{E} \subseteq \DiffR^n(F,R)$ such that 
		$$
		U = \big\lbrace f \in F \mid \delta(f) \in \pp \text{ for all } \delta \in \mathcal{E} \big\rbrace.
		$$
	\end{enumerate}
	\begin{proof}
		$(i)$ From \autoref{thm_Pimary_submod_as_Sol}$(ii.a)$, let $\mathcal{E}^\prime \subseteq \DiffR^n(F,R/\pp)$ be an  $(R \otimes_A R)$-submodule such that 
		$
		U = \Sol(\mathcal{E}^\prime).
		$
		Since $R$ is formally smooth and essentially of finite type over $A$ and $F$ is a finitely generated free $R$-module, \autoref{prop_exactness_diff_op} yields a canonical surjection 
		$$
		\Psi: \DiffR^n(F, R) \twoheadrightarrow \DiffR^n(F,R/\pp).
		$$
		So, the result follows by taking $\mathcal{E} = \Psi^{-1}\left(\mathcal{E}^\prime\right)$.
		
		$(ii)$ Follows identically to part $(i)$.
	\end{proof}
\end{corollary}

\section{A generalization of the Zariski-Nagata Theorem}
\label{sect_Zar_Nag}

In this section we provide a generalization of a celebrated theorem by Zariski and Nagata (see, e.g., \cite{ZARISKI}, \cite{NAGATA_LOCAL_RINGS}, \cite[Theorem 3.14]{EISEN_COMM}, \cite{RESUME_SYMB_HUNEKE_ET_AL}).
Here we relate symbolic powers to a new notion of differential powers, similarly to how it was done in \cite[\S 2.1]{RESUME_SYMB_HUNEKE_ET_AL}.
Although the outcomes of this section are simple consequences of \autoref{thm_Pimary_ideals_as_Sol}$(ii.b)$, they extend the recent results of \cite[Proposition 2.14]{RESUME_SYMB_HUNEKE_ET_AL}, \cite[Theorem 3.9]{ZAR_NAG_DESTEFANI_ET_AL} and \cite[Proposition 10.1]{brenner2018quantifying}.
Of particular interest is the fact that the new differential powers coincide with symbolic powers in non-smooth settings.

Throughout this section we use the following setup.

\begin{setup}
	\label{setup_symb_powers}
	Let $A$ be a Noetherian integral domain and $R$ be an $A$-algebra essentially of finite type over $A$ such that $A \subset R$.
\end{setup}

The $n$-th symbolic power of an ideal $I \subset R$ is defined as the ideal 
$$
I^{(n)}:= \bigcap_{\pp \in \Ass_R(R/I)} I^nR_\pp \cap R
$$
where $\pp$ runs through the associated primes of $R/I$.

\begin{remark}
	\label{rem_symb_power_smallest}
	It is well-known that the canonical map $R \rightarrow R_\pp$ gives a bijective correspondence between $\pp$-primary ideals in $R$ and $\pp R_\pp$-primary ideals in $R_\pp$ {\normalfont(}see, e.g., \cite[Theorem 4.1]{MATSUMURA}{\normalfont)}. 
	Therefore, for any $\pp \in \Spec(R)$, since $\pp^{(n)} = \pp^nR_\pp \cap R$, $\pp^{(n)}$ is the smallest $\pp$-primary ideal containing $\pp^n$.
\end{remark}

The corollary below describes symbolic powers of prime ideals as the solution sets of certain differential operators.

\begin{corollary}
	\label{cor_symb_powers_diff_ops}
	Assume \autoref{setup_symb_powers}. 
	Let $\pp \in \Spec(R)$ be a prime ideal in $R$ such that $\pp \cap A = 0$.
	Suppose that $\Quot(A) \hookrightarrow k(\pp)$ is a separable field extension, which holds whenever $\Quot(A)$ is perfect. 
	Let $N$ be a finitely generated torsion-free $R/\pp$-module.
	Then, for every $n \ge 1$, we have 
	$$
	\pp^{(n)} = \Sol\left(\DiffR^{n-1}(R, N)\right).
	$$
	\begin{proof}
		From \autoref{prop_sol_closed_submods_diff}$(iii)$ we obtain that $\Sol\left(\DiffR^{n-1}(R, N)\right)$ is a $\pp$-primary ideal containing $\pp^n$, then \autoref{rem_symb_power_smallest} implies that $\Sol\left(\DiffR^{n-1}(R, N)\right) \supseteq \pp^{(n)}$.
		Conversely, from \autoref{thm_Pimary_ideals_as_Sol}$(ii.b)$ we get an  $\left(R \otimes_A R\right)$-submodule $\mathcal{E} \subseteq \DiffR^{n-1}(R,N)$ such that $\pp^{(n)}=\Sol(\mathcal{E})$, and so it is clear that $\Sol\left(\DiffR^{n-1}(R, N)\right) \subseteq \pp^{(n)}$.
		So, the result follows.
	\end{proof}
\end{corollary}

To describe symbolic powers via differential operators, the following ideals have been studied in \cite{RESUME_SYMB_HUNEKE_ET_AL}.

\begin{definition}
	For $n \ge 1$ and $I \subset R$ an ideal in $R$, we set
	$$
	I^{{\langle n \rangle}_A} := \big\lbrace f \in R \mid \delta(f) \in I \text{ for all } \delta \in \DiffR^{n-1}(R,R)  \big\rbrace.
	$$
\end{definition}

Now, we define a new version of differential powers which seems to be better suited to describe symbolic powers, especially because it can be used in many interesting non-smooth situations.

\begin{definition}
	\label{def_new_diff_powers}
	For $n \ge 1$ and $I \subset R$ an ideal in $R$, we set
	$$
	I^{{\{n\}}_A} := \Sol\left(\DiffR^{n-1}\left(R,R/I\right)\right).
	$$
\end{definition}

The next theorem contains a generalization of the Zariski-Nagata Theorem, that relates the two above notions of differential powers with symbolic powers.

\begin{theorem}
	\label{thm_symb_powers_diff_powers}
	Assume \autoref{setup_symb_powers}. 
	Let $\pp \in \Spec(R)$ be a prime ideal in $R$ such that $\pp \cap A = 0$.
	Then, the following statements hold:
	\begin{enumerate}[(i)]
		\item $\pp^{(n)} \subseteq \pp^{{\{n\}}_A} \subseteq \pp^{{\langle n \rangle}_A}$.
		\item If $\Quot(A) \hookrightarrow k(\pp)$ is a separable field extension, which holds whenever $\Quot(A)$ is perfect, then
		$$
		\pp^{(n)} = \pp^{{\{n\}}_A}.
		$$ 
		\item If $R$ is formally smooth over $A$, then 
		$$
		\pp^{{\{n\}}_A} = \pp^{{\langle n \rangle}_A}.
		$$		
	\end{enumerate}
	\begin{proof}
		$(i)$
		The first inclusion $\pp^{(n)} \subseteq \pp^{{\{n\}}_A}$ is obtained from  \autoref{prop_sol_closed_submods_diff}$(iii)$ and \autoref{rem_symb_power_smallest}.
		For the second inclusion $\pp^{{\{n\}}_A} \subseteq \pp^{{\langle n \rangle}_A}$ we only need to note that the canonical surjection $R \twoheadrightarrow R/\pp$ induces a map 
		\begin{equation}
				\label{eq_map_diff_ops_reduction_prime}
				\DiffR^{n-1}(R, R) \rightarrow \DiffR^{n-1}(R, R/\pp)
		\end{equation}
		by post-composing any $\delta \in \DiffR^{n-1}(R,R)$ with the map $R \twoheadrightarrow R/\pp$; see \autoref{eq_Diff_functor}.	
		
		$(ii)$ It follows directly from \autoref{cor_symb_powers_diff_ops}.
		
		$(iii)$ If $R$ is formally smooth over $A$, then \autoref{prop_exactness_diff_op} implies that the map \autoref{eq_map_diff_ops_reduction_prime} above is surjective.
		So, the result is clear. 
	\end{proof}
\end{theorem}

\section{Some examples and computations}
\label{sect_examples}

In this section we include some simple examples. 
In order to make them illustrative, we shall try to make the computations self-contained  and without quoting the main results of the previous sections.

We start by describing ideals that are primary with respect to a maximal ideal that corresponds to a point in an affine space. 
This classical result is due to Gr\"obner (for more details, see, e.g., \cite{MARINARI_MORA}, \cite{MOURRAIN_DUALITY}, \cite[\S 10.2]{STURMFELS_SOLVING}).

\begin{example}
	\label{exam_grobner_duality}
	Let $\kk$ be a field and $R=\kk[x_1,\ldots,x_n]$. 
	Let $\alpha = \left(\alpha_1,\ldots,\alpha_n\right) \in \kk^n$ be a point and $\mm_\alpha=\left(x_1-\alpha_1, \ldots, x_n-\alpha_n\right) \subset R$ be the corresponding maximal ideal. 
	Then, for any $\mm_\alpha$-primary ideal $I \subset R$, there exists a finite number of differential operators $\delta_1, \ldots, \delta_m \in \Diff_{R/\kk}\left(R, R\right)$ such that 
	$$
	I = \big\{ f \in R \mid \delta_i(f)(\alpha_1,\ldots,\alpha_n) = 0 \text{ for all } 1 \le i \le m\big\}.
	$$
	\begin{proof}
		Suppose that $\mm_\alpha^{k+1} \subseteq I \subseteq \mm_\alpha$ for some $k \ge 0$. 
		Since $\kk\cong R/\mm_\alpha$, we consider $\kk$ as an $R$-module via the canonical homomorphism $R \twoheadrightarrow R/\mm_\alpha$.
		From \autoref{prop_represen_diff_opp} and the Hom-tensor adjunction we have the isomorphisms 
		$$
		\Diff_{R/\kk}^k\left(R, \kk\right) \cong \Hom_R\left(P_{R/\kk}^k, \kk\right) \cong \Hom_\kk\left(\kk \otimes_R P_{R/\kk}^k, \kk\right).
		$$
		We denote $R \otimes_\kk R = \kk[x_1, \ldots, x_n, \hat{x}_1, \ldots, \hat{x}_n]$ as a polynomial ring in $2n$ variables, where $x_i$ represents $x_i \otimes_\kk 1$ and $\hat{x}_i$ represents $1 \otimes_\kk x_i$.
		Thus, we get
		\begin{align*}
			\kk \otimes_R P_{R/\kk}^k &\cong \frac{\kk[x_1, \ldots, x_n, \hat{x}_1\ldots, \hat{x}_n]}{\Big(x_1-\alpha_1, \ldots, x_n-\alpha_n, {\left(x_1-\hat{x}_1, \ldots, x_n-\hat{x}_n\right)}^{k+1}\Big)} \\
			&\cong \frac{\kk[\hat{x}_1,\ldots,\hat{x}_n]}{{\left(\hat{x}_1-\alpha_1, \ldots, \hat{x}_n-\alpha_n\right)}^{k+1}}.
		\end{align*}
		Let $\hat{I} \subset \kk[\hat{x}_1, \ldots, \hat{x}_n]$ be the ideal obtained from $I$ by making the substitutions $x_i \mapsto \hat{x}_i$.
		Hence, we have the following inclusion 
		$$
		\Hom_\kk\left(\frac{\kk[\hat{x}_1,\ldots,\hat{x}_n]}{\hat{I}}, \kk\right) \hookrightarrow \Hom_\kk\left(\kk \otimes_R P_{R/\kk}^k, \kk\right) \cong \Diff_{R/\kk}^k\left(R, \kk\right),
		$$ 
		and one can see that its image yields an $\left(R \otimes_\kk R\right)$-submodule  $\mathcal{E}^\prime \subseteq \Diff_{R/\kk}^k(R, \kk)$ such that $I = \Sol\left(\mathcal{E}^\prime\right)$.
		From \cite[Example 4]{OBERST_NOETH_OPS} or \autoref{prop_exactness_diff_op}, the canonical map $R \twoheadrightarrow \kk$ gives the surjection
		$$
		\Psi : \Diff_{R/\kk}^k(R, R) \twoheadrightarrow \Diff_{R/\kk}^k(R, \kk).
		$$
		Take $\mathcal{E} =\Psi^{-1}\left(\mathcal{E}^\prime\right)$ and $\delta_1, \ldots, \delta_m$ a finite set of generators of $\mathcal{E}$.
		Therefore, we have
		$$
		I = \big\{ f \in R \mid \delta_i(f) \in \mm_\alpha  \text{ for all } 1 \le i \le m\big\},
		$$
		and so the result follows.
	\end{proof}
\end{example}

The purpose of the next example is twofold.
It shows that without any assumption of smoothness: the symbolic powers may not be equal to the differential powers studied in \cite{RESUME_SYMB_HUNEKE_ET_AL} but may still coincide with the new notion of differential powers (\autoref{def_new_diff_powers}), and there may exist primary ideals which cannot be described as in \autoref{cor_prim_ideals_sol_Diff_R_R}.

\begin{example}
	\label{exam_cubic}
	Let $\CC$ be the field of complex numbers, $R = \CC[x,y,z]/(x^3+y^3+z^3)$ and $\mm=(x,y,z)R$ be the graded irrelevant ideal of $R$.
	Then, for all $n \ge 2$, the following statements hold:
	\begin{enumerate}[(i)]
		\item $\mm^{{\langle n \rangle}_\CC} = \mm$.
		\item $\mm^n$ cannot be described as in  \autoref{cor_prim_ideals_sol_Diff_R_R}. 
		\item $\mm^{{\{ n \}}_\CC} = \mm^n$.
	\end{enumerate}
	\begin{proof}
		$(i), (ii)$
		Note that $R$ is a standard graded $\CC$-algebra.
		From \cite{BGG_NON_NOETHERIAN_DIFF} we know that $\Diff_{R/\CC}(R, R)$ is not a finitely generated $\CC$-algebra and also not Noetherian, that $\Diff_{R/\CC}(R, R)$ is graded where an operator $\delta \in \Diff_{R/\CC}(R, R)$ is homogeneous of degree $k$ if it satisfies the condition 
		$$
		\delta\left({\left[R\right]}_i\right) \subseteq {\left[R\right]}_{i+k}
		$$
		for all $i \in \ZZ$, and that with this grading ${\left[\Diff_{R/\CC}(R, R)\right]}_k=0$ for all $k < 0$.
		
		The above remarks imply that for any $\delta \in \Diff_{R/\CC}(R, R)$, we have that $\delta(\mm) \subseteq \mm$.
		So, for any $(R\otimes_\CC R)$-submodule $\mathcal{E} \subseteq \Diff_{R/\CC}(R, R)$ we obtain
		$$
		\big\{ f \in R \mid \delta(f) \in \mm \text{ for all } \delta \in \mathcal{E} \big\} \supseteq \mm,
		$$
		and this implies the statement of part $(ii)$ and that $\mm^{{\langle n\rangle}_\CC} \supseteq \mm$.
		Since $1 \in R \cong \Diff_{R/\CC}^{0}(R, R)  \subseteq \Diff_{R/\CC}^{n-1}(R, R)$, it follows that $\mm^{{\langle n\rangle}_\CC}=\mm$.
		
		$(iii)$ Here the argument is completely similar to \autoref{exam_grobner_duality}.
		Again, we see $\CC$ as an $R$-module via the canonical map $R \twoheadrightarrow R/\mm \cong \CC$.
		For any $n \ge 1$, \autoref{prop_represen_diff_opp} and the Hom-tensor adjunction yield the isomorphisms
		\begin{align*}
			\Diff_{R/\CC}^{n-1}(R, \CC) &\cong \Hom_R\left(P_{R/\CC}^{n-1}, \CC\right) \\
			&\cong \Hom_\CC\left(\CC \otimes_R P_{R/\CC}^{n-1}, \CC\right)\\
			&\cong 	\Hom_\CC\left(\frac{\CC[\hat{x}, \hat{y}, \hat{z}]}{\big(\hat{x}^3+\hat{y}^3+\hat{z}^3, {(\hat{x}, \hat{y}, \hat{z})}^{n}\big)}, \CC\right),		
		\end{align*}
		where we are using the notation 
		$$
		R \otimes_\CC R \cong \frac{\CC[x,y,z, \hat{x}, \hat{y}, \hat{z}]}{\left(x^3+y^3+z^3, \hat{x}^3+\hat{y}^3+\hat{z}^3\right)}.
		$$
		Therefore, it is clear that $\mm^{{\{ n \}}_\CC} =\Sol\left(\Diff_{R/\CC}^{n-1}(R, \CC)\right) = \mm^n$.
	\end{proof}
\end{example}

Finally, the last example shows that the condition of separability in \autoref{thm_symb_powers_diff_powers}$(ii)$ cannot be avoided (also, see \cite[Example 3.8]{ZAR_NAG_DESTEFANI_ET_AL}).

\begin{example}
	Let $p \in \NN$ be a prime, $\kk=\bFF_p(t)$, $R = \kk[x]$ and $\pp = (x^p - t)$.
	Then, we have that $\pp^{{\{2\}}_\kk} = \pp \neq \pp^2 = \pp^{(2)}$.
	\begin{proof}
		Set $R/\pp \cong \kk(u)$ for some $u \in \overline{\kk}$ (in an algebraic closure of $\kk$) that satisfies the equation $u^p-t=0$.
		Again, we use the notation $R \otimes_\kk R = \kk[x,\hat{x}]$ and obtain the isomorphisms 
		\begin{align*}
		\Diff_{R/\kk}^{1}(R, \kk(u)) &\cong \Hom_R\left(P_{R/\kk}^{1}, \kk(u)\right)\\ 
		&\cong \Hom_{\kk(u)}\left(\kk(u) \otimes_R P_{R/\kk}^{1}, \kk(u)\right)\\
		&\cong 	\Hom_{\kk(u)}\left(\frac{\kk(u)[\hat{x}]}{{(\hat{x}-u)}^{2}}, \kk(u)\right)		
		\end{align*}
		(see \autoref{exam_grobner_duality} and \autoref{exam_cubic}).
		Since $\left(\hat{x}^p - t\right)={\left(\hat{x}-u\right)}^p \subseteq {\left(\hat{x}-u\right)}^2$, it clearly follows that $\pp^{{\{2\}}_\kk} = \Sol\left(\Diff_{R/\kk}^{1}(R, \kk(u))\right) \supseteq \pp$.
		Since the canonical map $R \twoheadrightarrow R/\pp$ belongs to $\Diff_{R/\CC}^{0}(R, R/\pp)  \subseteq \Diff_{R/\CC}^{1}(R, R/\pp)$, we get that $\pp^{{\{2\}}_\kk} = \pp$.
		Finally, we have that $\pp^{(2)}=\pp^2$, because $\pp$ is maximal.
	\end{proof} 
\end{example}

\section*{Acknowledgments}

The author is grateful to Bernd Sturmfels for suggesting the initial problem that led to the preparation of this paper, for many helpful discussions, and for his contagious enthusiasm.
The author thanks {\L}ukasz Grabowski and Viktor Levandovskyy for useful conversations. 
The author thanks Greg Brumfiel for a helpful correspondence.
The author is grateful to the referee for valuable comments and suggestions for the improvement of this work.

\bibliographystyle{elsarticle-num} 

\end{document}